\documentclass[10pt]{amsart}
\usepackage{amssymb,latexsym}

\usepackage{ifthen}
\usepackage{graphicx,color}
\usepackage{epic,eepic}

\usepackage{microtype}

\newcommand{\light}[1]{{#1}}
\newcommand{\dark}[1]{{#1}}
\newcommand{\dgreen}[1]{{#1}}

\newcounter{colorversion}

\ifthenelse{\thecolorversion=1} {
\usepackage{color}
\definecolor{light}{rgb}{1,0.1,0.1}        
\definecolor{dark}{cmyk}{1,0.4,0,0}  
\definecolor{dgreen}{rgb}{0.00,0.59,0.00}  
\definecolor{bulletcolor}{cmyk}{1,0.3,1,0}  
\renewcommand{\light}[1]{\color{light}{#1}\color{black}}
\renewcommand{\dark}[1]{\color{dark}{#1}\color{black}}
\renewcommand{\dgreen}[1]{\color{dgreen}{#1}\color{black}}

}{}

\numberwithin{equation}{section}

\newtheorem{thm}{Theorem}[section]
\newtheorem{prop}[thm]{Proposition}

\theoremstyle{definition}

\newtheorem{example}[thm]{Example}

\def\endproof{\hfill$\square$\medskip}

\def\ZZ{\mathbb{Z}}
\def\CC{\mathbb{C}}

\def\lg{\mathfrak{g}}
\def\nn{\mathfrak{n}}
\def\hh{\mathfrak{h}}

\def\gg{\mathbf{g}}

\def\ii{\mathbf{i}}
\def\vv{\mathbf{v}}

\def\sp{\mathfrak{sp}}
\def\so{\mathfrak{so}}
\def\soo{\mathfrak{so}_{2n+1}}
\def\soe{\mathfrak{so}_{2n}}
\def\SL{\mathrm{SL}}

\def\l{\ell}
\def\wnot{w_\circ}

\newcommand{\barred}[1]{{\,\overline{#1}\,}}

\newcommand{\bmcase}{\left\{ \begin{array}{ll}}
\newcommand{\emcase}{\end{array}\right.}

\newcommand{\beal}{\begin{eqnarray}\begin{array}{l} }
\newcommand{\eea}{\end{array}\end{eqnarray}}

\renewcommand{\eqref}[1]{{\rm (\ref{#1})}}

\newcommand{\mat}[4]{\left(\!\!\begin{array}{cc}
#1 & #2 \\
#3 & #4 \\
\end{array}\!\!\right)}

\begin{document}

\title{Combinatorial Expressions for $F$-polynomials in Classical Types}

\author{Shih-Wei Yang}
\address{\noindent Department of Mathematics, Northeastern University,
 Boston, MA 02115}
\email{yang.s@neu.edu}

\begin{abstract}
We give combinatorial formulas for $F$-polynomials in cluster algebras of classical types in terms of the weighted paths in certain directed graphs. As a consequence we prove the positivity of $F$-polynomials in cluster algebras of classical types.
\end{abstract}

\subjclass[2010]{Primary
 13F60. 
       }


\thanks{The author is supported by A.~Zelevinsky's NSF grant \#DMS-0801187.}

\maketitle

\tableofcontents

\section{Introduction and main results}
\label{sec:introduction}

$F$-polynomials are a important family of polynomials in the theory of cluster algebras. As shown in~\cite{ca4}, $F$-polynomials and $\gg$-vectors provide all the information needed to express \emph{any} cluster variable in a cluster algebra with \emph{arbitrary} coefficients in terms of an initial cluster. It was shown in~\cite{yangzel} that in a cluster algebra of finite type with an acyclic initial seed, the $F$-polynomials are given by a certain set of generalized principal minors. Generalized minors, first introduced in \cite{fz-double} for the study of total positivity in a simply connected semisimple complex algebraic group $G$, are a special family of regular functions $\Delta_{\gamma, \delta}$ on $G$. These functions are suitably normalized matrix coefficients corresponding to pairs of extremal weights $(\gamma, \delta)$ in some fundamental representation of $G$. we call a generalized minor $\Delta_{\gamma, \delta}$ \emph{principal} if $\gamma = \delta$.

Let $G$ be a simply connected semisimple complex Lie group with rank $n$. For $i \in [1, n] = \{1, \ldots, n\}$, let $\varphi_i : \SL_2 \rightarrow G$ denote the canonical embedding corresponding to the simple root $\alpha_i$. For $i \in [1, n]$ and $t \in \CC$, we write
\beal\label{eq:x,y:noef}
x_i (t) = \varphi_i \mat{1}{t}{0}{1}\,, \quad x_{\barred{i}} (t) = \varphi_i \mat{1}{0}{t}{1}\,.
\eea
Let $W$ be the Weyl group of $G$. Recall that~$W$ is a finite Coxeter group generated by the simple reflections $s_i$ for $i \in [1,n]$, and a Coxeter element $c$ is a element in $W$ such that $c = s_{i_1} \cdots s_{i_n}$ for some permutation $(i_1, \dots, i_n)$ of the index set~$[1,n]$. We also denote the longest element of~$W$ by $w_\circ$.

It was shown in~\cite{yangzel} that in a cluster algebra of an arbitrary finite type with arbitrary acyclic initial seed (depending on the choose of a Coxeter element $c = s_{i_1} \cdots s_{i_n}$), the $F$-polynomials can be parametrized by a special set of extremal weights $\{c^m \omega_k : k \in [1,n],\, 0 \leq m \leq h(k;c)\}$ where $h(k;c)$ is the smallest positive integer such that $c^{h(k;c)} \omega_k = w_\circ (\omega_k)$.
According to~\cite[Theorem~1.12]{yangzel}, the $F$-polynomials are given by
\begin{equation}
\label{eq:F-poly}
F_{c^m \omega_k}(t_1, \dots, t_n) =
\Delta_{c^m \omega_k, c^m \omega_k}
(x_{\barred{i_i}}(1) \cdots x_{\barred{i_n}}(1) x_{i_n}(t_{i_n}) \cdots x_{i_1}(t_{i_1})).
\end{equation}
Our main result is the following theorem.
\begin{thm}
\label{th:pF-poly}
In the cluster algebra of classical type with an arbitrary acyclic initial seed, explicit combinatorial expressions for the $F$-polynomials are given. The descriptions for the types $\mathsf{A_n}$, $\mathsf{D_n}$, $\mathsf{B_n}$ and $\mathsf{C_n}$ are given in Propositions~\ref{pr:An-minors}, ~\ref{pr:genminorsD}, ~\ref{pr:genminorsB} and~\ref{pr:genminorsC} respectively. Furthermore, in all these cases, the coefficients of the $F$-polynomials are manifestly positive.
\end{thm}

There are other formulas for the $F$-polynomials and proofs for the positivity conjecture in the literature. In particular, Fomin and Zelevinsky's work in~\cite{yga} together with~\cite{ca4} gave explicit formulas and proved the positivity for the $F$-polynomials in classical types for a bipartite initial cluster; Musiker, Schiffler and Williams's work in~\cite{msw} deals with cluster algebras from surfaces; the results in~\cite{tran1, tran2} by Tran have the same generality as this current work. Our answer is given in very different terms and obtained by totally different methods.

The proof of Theorem~\ref{th:pF-poly} is based on combinatorial formulas for generalized minors in the classical types of the form given in~(\ref{eq:F-poly}). In the type $\mathsf{A_n}$ case, the generalized minors specialize to the ordinary minors and our combinatorial formula is a well-known result due to Lindstr\"{o}m (see~\cite{lindstrom}, \cite{GV}, \cite{GV2}, \cite{fz-double} and \cite{fz-intel}). For the convenience of the reader, we will recall the type $\mathsf{A_n}$ theory in this note.
For the type $\mathsf{B_n}$, we will construct two weighted directed graphs $\Gamma (\mathsf{B_n}, c)$ and $\Gamma_{\mathrm{S}} (\mathsf{B_n}, c)$, while for the types $\mathsf{C_n}$ and $\mathsf{D_n}$, we only need one directed graph for each type, $\Gamma (\mathsf{C_n}, c)$ and $\Gamma (\mathsf{D_n}, c)$ respectively. The formulas are given in terms of the weighted paths in the corresponding directed graphs. All the proofs of the results in this section will be given in Section~\ref{sec:proofs}.

\textbf{Type $\mathsf{A_n}$ :}
Let $E_{i,j}$ denote the $(n+1) \times (n+1)$ matrix whose
$(i,j)$-entry is equal to~1 while all other entries are~0, and
let $\mathrm{Id} \in G$ denote the identity matrix.
For $i=1,\dots,n$, let

\begin{eqnarray}\begin{array}{l}
\label{eq:jacobi}
x_i(t)=\mathrm{Id}+tE_{i,i+1}=
\left(\begin{array}{cccccc}
1      & \cdots & 0      & 0      & \cdots & 0      \\
\cdots & \cdots & \cdots & \cdots & \cdots & \cdots \\
0      & \cdots & 1      & t      & \cdots & 0      \\
0      & \cdots & 0      & 1      & \cdots & 0      \\
\cdots & \cdots & \cdots & \cdots & \cdots & \cdots \\
0      & \cdots & 0      & 0      & \cdots & 1 \\
\end{array}\right)
\end{array}\end{eqnarray}
and
\begin{eqnarray}\begin{array}{l}
x_{\barred{i}}(t)=\mathrm{Id}+tE_{i+1,i}=
\left(\begin{array}{cccccc}
1      & \cdots & 0      & 0      & \cdots & 0      \\
\cdots & \cdots & \cdots & \cdots & \cdots & \cdots \\
0      & \cdots & 1      & 0      & \cdots & 0      \\
0      & \cdots & t      & 1      & \cdots & 0      \\
\cdots & \cdots & \cdots & \cdots & \cdots & \cdots \\
0      & \cdots & 0      & 0      & \cdots & 1 \\
\end{array}\right)\,.
\end{array}\end{eqnarray}

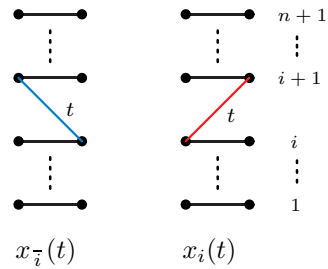
\begin{figure}[h]
\setlength{\unitlength}{1.2pt}
\begin{center}
\begin{picture}(50,80)(0,-15)
\thicklines
\put(5,0){\line(1,0){20}} \put(5,20){\line(1,0){20}} \put(5,40){\line(1,0){20}} \put(5,60){\line(1,0){20}}
  \put(5,0){\circle*{2.5}}
  \put(5,20){\circle*{2.5}}
  \put(5,40){\circle*{2.5}}
  \put(5,60){\circle*{2.5}}
  \put(25,0){\circle*{2.5}}
  \put(25,20){\circle*{2.5}}
  \put(25,40){\circle*{2.5}}
  \put(25,60){\circle*{2.5}}
  \dottedline{3}(15,5)(15,15)
  \dottedline{3}(15,45)(15,55)
\dark{ \put(5,40){\line(1,-1){20}} }
  \put(20,28){\footnotesize$t$}
  \put(4,-17){$x_{\barred{i}}(t)$}
\end{picture}
\begin{picture}(50,80)(0,-15)
\thicklines
\put(5,0){\line(1,0){20}} \put(5,20){\line(1,0){20}} \put(5,40){\line(1,0){20}} \put(5,60){\line(1,0){20}}
  \put(5,0){\circle*{2.5}}
  \put(5,20){\circle*{2.5}}
  \put(5,40){\circle*{2.5}}
  \put(5,60){\circle*{2.5}}
  \put(25,0){\circle*{2.5}}
  \put(25,20){\circle*{2.5}}
  \put(25,40){\circle*{2.5}}
  \put(25,60){\circle*{2.5}}
  \dottedline{3}(15,5)(15,15)
  \dottedline{3}(15,45)(15,55)
  \put(38,-2){\scriptsize${1}$}
  \dottedline{3}(40,7)(40,14)
  \put(38,18){\scriptsize${i}$}
  \put(34,38){\scriptsize${i+1}$}
  \dottedline{3}(40,47)(40,53)
  \put(34,58){\scriptsize${n+1}$}
\light{ \put(5,20){\line(1,1){20}} }
  \put(18,26){\footnotesize$t$}
  \put(4,-17){$x_{i}(t)$}
\end{picture}
\end{center}
\caption{``Elementary chips'' of type $\mathsf{A_n}$\,.} \label{fig:Achips}
\end{figure}

For any $i \in [1, n] \cup [\barred{1}, \barred{n}]$, where $[\barred{1}, \barred{n}] = \{\barred{1}, \ldots, \barred{n}\}$, we construct an ``elementary chip'' corresponding to $x_{i} (t)$ to be a weighted directed graph of one of the kinds shown in Figure~\ref{fig:Achips}.

Note that in each chip, the horizontal levels are labeled by $1, \ldots, n$ starting from the bottom.
The chip corresponding to $x_{i}(t)$ or $x_{\barred{i}}(t)$ has a diagonal edge connecting the horizontal levels $i$ and $i+1$ with weight $t$. All other (unlabeled) edges have weight~1 and all edges are presumed to be oriented from right to left.
The directed graph $\Gamma (\mathsf{A_n}, c)$ associated with $c = s_{i_1} \cdots s_{i_n}$ is constructed as a concatenation of elementary chips $x_{\barred{i_i}}(1), \ldots, x_{\barred{i_n}}(1), x_{i_n}(t_{i_n}), \ldots, x_{i_1}(t_{i_1})$ (in this order). We number the $n+1$ sources and $n+1$ sinks of the graph $\Gamma (\mathsf{A_n}, c)$ bottom-to-top, and define the weight of a path in~$\Gamma (\mathsf{A_n}, c)$ to be the product of the weights of all edges in the path. We also define the weight of a family of paths to be the product of the weights of all paths in the family.

The Weyl group of type $\mathsf{A_n}$ is identified with the symmetric group $S_{n+1}$ and it acts on the index set $[1, n+1]$ as permutations. The simple reflections are $s_i = (i, i+1)$ for $i \in [1,n]$. Then the $F$-polynomials in type $\mathsf{A_n}$ are computed as follows:

\begin{prop}
\label{pr:An-minors}
The $F$-polynomial 
$F_{c^m \omega_k}(t_1, \dots, t_n)$
equals the sum of weights of all collections of vertex-disjoint paths in $\Gamma (\mathsf{A_n}, c)$ with the sources and sinks labeled by $c^m \cdot [1,k]$.
\end{prop}

\begin{example}
Type $\mathsf{A_3}$: Let $c = s_1s_3s_2 = (1, 2, 4, 3)$, then $c \cdot[1, 2]=\{2, 4\}$, hence
$$F_{c \omega_2}(t_1, t_2, t_3) = 1 + t_1 + t_3 + t_1t_3 + t_1 t_2t_3\,.$$
In Figure~\ref{fig:A3ex}, we give all families of vertex-disjoint paths in $\Gamma (\mathsf{A_3}, s_1s_3s_2)$ with the sources and sinks labeled by $\{2, 4\}$ and each family of paths is depicted by thick lines.

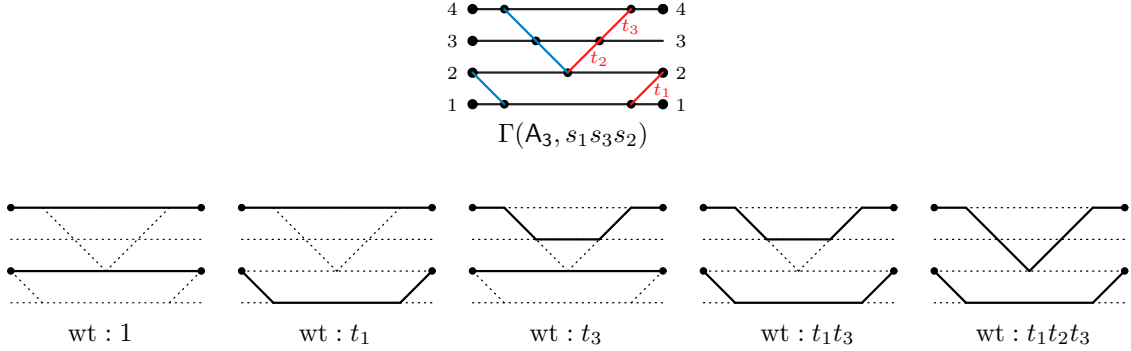
\begin{figure}[h!]
\setlength{\unitlength}{1.2pt}
\begin{center}
\begin{picture}(70,40)(0,-10)
\thicklines
\drawline[90](0,0)(60,0)
\drawline[90](0,10)(60,10)
\drawline[90](0,20)(60,20)
\drawline[90](0,30)(60,30)
  \put(0,0){\circle*{2.5}}
  \put(0,10){\circle*{2.5}}
  \put(0,20){\circle*{2.5}}
  \put(0,30){\circle*{2.5}}
  \put(60,0){\circle*{2.5}}
  \put(60,10){\circle*{2.5}}
  \put(60,30){\circle*{2.5}}
  \put(60,30){\circle*{2.5}}
\put(50,0){\circle*{2}}
\put(50,30){\circle*{2}}
\put(40,20){\circle*{2}}
\put(30,10){\circle*{2}}
\put(20,20){\circle*{2}}
\put(10,0){\circle*{2}}
\put(10,30){\circle*{2}}
  \put(-8,-2){\scriptsize$1$}
  \put(-8,8){\scriptsize$2$}
  \put(-8,18){\scriptsize$3$}
  \put(-8,28){\scriptsize$4$}
  \put(64,-2){\scriptsize$1$}
  \put(64,8){\scriptsize$2$}
  \put(64,18){\scriptsize$3$}
  \put(64,28){\scriptsize$4$}
\light{
\drawline[90](60,10)(50,0)
\drawline[90](50,30)(40,20)
\drawline[90](40,20)(30,10)
  \put(57,3){\scriptsize$t_1$}
  \put(37,13){\scriptsize$t_2$}
  \put(47,23){\scriptsize$t_3$}
}
\dark{
\drawline[90](30,10)(20,20)
\drawline[90](20,20)(10,30)
\drawline[90](10,0)(0,10)
}
\put(8,-12){$\Gamma (\mathsf{A_3}, s_1s_3s_2)$}
\end{picture}
\end{center}

\begin{center}
\begin{picture}(70,60)(0,-10)
\thinlines
\dottedline{2}(0,0)(60,0)
\dottedline{2}(0,10)(60,10)
\dottedline{2}(0,20)(60,20)
\dottedline{2}(0,30)(60,30)
  \put(0,30){\circle*{2}}
  \put(0,10){\circle*{2}}
  \put(60,30){\circle*{2}}
  \put(60,10){\circle*{2}}
\dottedline{2}(60,10)(50,0)
\dottedline{2}(50,30)(30,10)(10,30)
\dottedline{2}(10,0)(0,10)
\thicklines
\drawline[90](0,30)(60,30)
\drawline[90](0,10)(60,10)
\put(18,-12){$\mathrm{wt}: 1$}
\end{picture}
\begin{picture}(70,60)(0,-10)
\thinlines
\dottedline{2}(0,0)(60,0)
\dottedline{2}(0,10)(60,10)
\dottedline{2}(0,20)(60,20)
\dottedline{2}(0,30)(60,30)
  \put(0,30){\circle*{2}}
  \put(0,10){\circle*{2}}
  \put(60,30){\circle*{2}}
  \put(60,10){\circle*{2}}
\dottedline{2}(60,10)(50,0)
\dottedline{2}(50,30)(30,10)(10,30)
\dottedline{2}(10,0)(0,10)
\thicklines
\drawline[90](0,30)(60,30)
\drawline[90](0,10)(10,0)(50,0)(60,10)
\put(18,-12){$\mathrm{wt}: t_1$}
\end{picture}
\begin{picture}(70,60)(0,-10)
\thinlines
\dottedline{2}(0,0)(60,0)
\dottedline{2}(0,10)(60,10)
\dottedline{2}(0,20)(60,20)
\dottedline{2}(0,30)(60,30)
  \put(0,30){\circle*{2}}
  \put(0,10){\circle*{2}}
  \put(60,30){\circle*{2}}
  \put(60,10){\circle*{2}}
\dottedline{2}(60,10)(50,0)
\dottedline{2}(50,30)(30,10)(10,30)
\dottedline{2}(10,0)(0,10)
\thicklines
\drawline[90](0,30)(10,30)(20,20)(40,20)(50,30)(60,30)
\drawline[90](0,10)(60,10)
\put(18,-12){$\mathrm{wt}: t_3$}
\end{picture}
\begin{picture}(70,60)(0,-10)
\thinlines
\dottedline{2}(0,0)(60,0)
\dottedline{2}(0,10)(60,10)
\dottedline{2}(0,20)(60,20)
\dottedline{2}(0,30)(60,30)
  \put(0,30){\circle*{2}}
  \put(0,10){\circle*{2}}
  \put(60,30){\circle*{2}}
  \put(60,10){\circle*{2}}
\dottedline{2}(60,10)(50,0)
\dottedline{2}(50,30)(30,10)(10,30)
\dottedline{2}(10,0)(0,10)
\thicklines
\drawline[90](0,30)(10,30)(20,20)(40,20)(50,30)(60,30)
\drawline[90](0,10)(10,0)(50,0)(60,10)
\put(18,-12){$\mathrm{wt}: t_1t_3$}
\end{picture}
\begin{picture}(70,60)(0,-10)
\thinlines
\dottedline{2}(0,0)(60,0)
\dottedline{2}(0,10)(60,10)
\dottedline{2}(0,20)(60,20)
\dottedline{2}(0,30)(60,30)
  \put(0,30){\circle*{2}}
  \put(0,10){\circle*{2}}
  \put(60,30){\circle*{2}}
  \put(60,10){\circle*{2}}
\dottedline{2}(60,10)(50,0)
\dottedline{2}(50,30)(30,10)(10,30)
\dottedline{2}(10,0)(0,10)
\thicklines
\drawline[90](0,30)(10,30)(30,10)(50,30)(60,30)
\drawline[90](0,10)(10,0)(50,0)(60,10)
\put(13,-12){$\mathrm{wt}: t_1t_2t_3$}
\end{picture}
\end{center}
\caption{$\Gamma(\mathsf{A_3}, s_1s_3s_2)$, $F_{c \omega_2}(t_1, t_2, t_3)$.} \label{fig:A3ex}
\end{figure}
\end{example}

\begin{figure}[h]
\setlength{\unitlength}{1.2pt}
\begin{center}
\begin{picture}(63,120)(0,-15)
\thicklines
\put(18,0){\line(1,0){20}} \put(18,20){\line(1,0){20}} \put(18,40){\line(1,0){20}} \put(18,60){\line(1,0){20}}
\put(18,80){\line(1,0){20}} \put(18,100){\line(1,0){20}}
  \put(18,0){\circle*{2.5}}
  \put(18,20){\circle*{2.5}}
  \put(18,40){\circle*{2.5}}
  \put(18,60){\circle*{2.5}}
  \put(18,80){\circle*{2.5}}
  \put(18,100){\circle*{2.5}}
  \put(38,0){\circle*{2.5}}
  \put(38,20){\circle*{2.5}}
  \put(38,40){\circle*{2.5}}
  \put(38,60){\circle*{2.5}}
  \put(38,80){\circle*{2.5}}
  \put(38,100){\circle*{2.5}}
  \put(4,-2){\scriptsize${1}$}
  \put(4,18){\scriptsize${i}$}
  \put(0,38){\scriptsize${i+1}$}
  \put(0,58){\scriptsize${\overline{i+1}}$}
  \put(4,78){\scriptsize${\overline{i}}$}
  \put(4,98){\scriptsize${\overline{1}}$}
  \dottedline{3}(6,6)(6,13)
  \dottedline{3}(6,46)(6,53)
  \dottedline{3}(6,86)(6,93)
  \dottedline{3}(28,5)(28,15)
  \dottedline{3}(28,45)(28,55)
  \dottedline{3}(28,85)(28,95)
\dark{ \put(18,40){\line(1,-1){20}} \put(18,80){\line(1,-1){20}} }
  \put(22,26){\footnotesize$t$}
  \put(22,66){\footnotesize$t$}
  \put(17,-17){$x_{\barred{i}} (t)$}
\end{picture}
\begin{picture}(50,120)(0,-15)
\thicklines
\put(5,0){\line(1,0){20}} \put(5,20){\line(1,0){20}} \put(5,40){\line(1,0){20}} \put(5,60){\line(1,0){20}}
\put(5,80){\line(1,0){20}} \put(5,100){\line(1,0){20}}
  \put(5,0){\circle*{2.5}}
  \put(5,20){\circle*{2.5}}
  \put(5,40){\circle*{2.5}}
  \put(5,60){\circle*{2.5}}
  \put(5,80){\circle*{2.5}}
  \put(5,100){\circle*{2.5}}
  \put(25,0){\circle*{2.5}}
  \put(25,20){\circle*{2.5}}
  \put(25,40){\circle*{2.5}}
  \put(25,60){\circle*{2.5}}
  \put(25,80){\circle*{2.5}}
  \put(25,100){\circle*{2.5}}
  \dottedline{3}(15,5)(15,15)
  \dottedline{3}(15,45)(15,55)
  \dottedline{3}(15,85)(15,95)
\light{ \put(5,20){\line(1,1){20}} \put(5,60){\line(1,1){20}} }
  \put(18,26){\footnotesize$t$}
  \put(18,66){\footnotesize$t$}
  \put(4,-17){$x_i (t)$}
\end{picture}
\begin{picture}(50,120)(0,-15)
\thicklines
\put(5,0){\line(1,0){20}} \put(5,20){\line(1,0){20}} \put(5,40){\line(1,0){20}} \put(5,60){\line(1,0){20}}
\put(5,80){\line(1,0){20}} \put(5,100){\line(1,0){20}}
  \put(5,0){\circle*{2.5}}
  \put(5,20){\circle*{2.5}}
  \put(5,40){\circle*{2.5}}
  \put(5,60){\circle*{2.5}}
  \put(5,80){\circle*{2.5}}
  \put(5,100){\circle*{2.5}}
  \put(25,0){\circle*{2.5}}
  \put(25,20){\circle*{2.5}}
  \put(25,40){\circle*{2.5}}
  \put(25,60){\circle*{2.5}}
  \put(25,80){\circle*{2.5}}
  \put(25,100){\circle*{2.5}}
  \dottedline{3}(15,5)(15,15)
  \dottedline{3}(15,85)(15,95)
\dark{ \put(5,60){\line(1,-2){20}} \put(5,80){\line(1,-2){20}} }
  \put(15,66){\footnotesize$t$}
  \put(14,27){\footnotesize$t$}
  \put(4,-17){$x_{\barred{n}} (t)$}
\end{picture}
\begin{picture}(50,120)(0,-15)
\thicklines
\put(5,0){\line(1,0){20}} \put(5,20){\line(1,0){20}} \put(5,40){\line(1,0){20}} \put(5,60){\line(1,0){20}}
\put(5,80){\line(1,0){20}} \put(5,100){\line(1,0){20}}
  \put(5,0){\circle*{2.5}}
  \put(5,20){\circle*{2.5}}
  \put(5,40){\circle*{2.5}}
  \put(5,60){\circle*{2.5}}
  \put(5,80){\circle*{2.5}}
  \put(5,100){\circle*{2.5}}
  \put(25,0){\circle*{2.5}}
  \put(25,20){\circle*{2.5}}
  \put(25,40){\circle*{2.5}}
  \put(25,60){\circle*{2.5}}
  \put(25,80){\circle*{2.5}}
  \put(25,100){\circle*{2.5}}
  \put(35,-2){\scriptsize${1}$}
  \put(30,18){\scriptsize${n-1}$}
  \put(35,38){\scriptsize${n}$}
  \put(35,58){\scriptsize${\overline{n}}$}
  \put(30,78){\scriptsize${\overline{n-1}}$}
  \put(35,98){\scriptsize${\overline{1}}$}
  \dottedline{3}(15,5)(15,15)
  \dottedline{3}(15,85)(15,95)
  \dottedline{3}(36,8)(36,15)
  \dottedline{3}(36,88)(36,95)
\light{ \put(5,20){\line(1,2){20}} \put(5,40){\line(1,2){20}} }
  \put(13,66){\footnotesize$t$}
  \put(14,27){\footnotesize$t$}
  \put(4,-17){$x_n (t)$}
\end{picture}
\end{center}
\caption{Elementary chips of type $\mathsf{D_n}$ ($i = 1, \ldots, n-1$).} \label{fig:dnchipsnsp}
\end{figure}
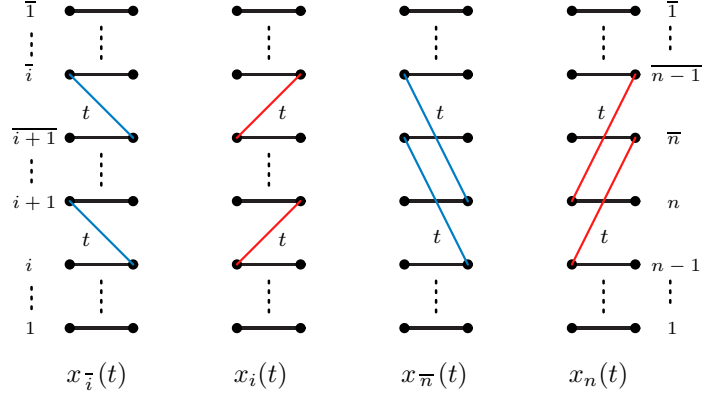

\textbf{Type $\mathsf{D_n}$} ($n \geq 4$) :
We use the standard numbering of simple roots as in~\cite{bourbaki}. For each $i \in [1, n] \cup [\barred{1}, \barred{n}]$, the elementary chip corresponding to $x_{i}(t)$ is shown in Figure~\ref{fig:dnchipsnsp}.
In each chip, the vertices consist of all the endpoints of the horizontal edges and all of the edges are oriented from right to left. We number the horizontal levels from bottom to top in the order  $1, \ldots, n, \barred{n}, \ldots, \barred{1}$. The numbering of the horizontal levels for the first (resp., last) two chips in Figure~\ref{fig:dnchipsnsp} is shown on the left (resp., right) of the figure. The two diagonal edges in each chip have weight $t$, all other unlabeled edges have weight 1.

The directed graph $\Gamma (\mathsf{D_n}, c)$ associated with $c = s_{i_1} \cdots s_{i_n}$ is constructed as a concatenation of elementary chips $x_{\barred{i_i}}(1), \ldots, x_{\barred{i_n}}(1), x_{i_n}(t_{i_n}), \ldots, x_{i_1}(t_{i_1})$ (in this order). We number the $2n$ sources and the $2n$ sinks of the graph $\Gamma (\mathsf{D_n}, c)$ bottom-to-top in the order $1, \ldots, n, \barred{n}, \ldots, \barred{1}$.

Note that in the chips corresponding to $x_{\barred{n}}(t)$ and $x_n(t)$, the intersections of the diagonal edges and the horizontal edges in the middle of each horizontal edge are \emph{not} vertices.  We call a family of paths \emph{bundled} if within each elementary chip, either both of the diagonal edges belong to the family of paths, or neither belong to the family of paths.
The weight of a family of paths is defined in the same way as in the type $\mathsf{A_n}$.

The Weyl group of type $\mathsf{D_n}$ acts on the index set $[1, n]$ as permutations with even number of ``bar'' changes. 
When written as permutations on $[1, n] \cup [\barred{1}, \barred{n}]$, the simple reflections are $s_i = (i, i+1)(\overline{i+1}, \overline{i})$ for $i = 1, \ldots, n-1$, and $s_n = (n-1, \overline{n})(n, \overline{n-1})$. Then the $F$-polynomials in type $\mathsf{D_n}$ are computed as follows:

\begin{prop}
\label{pr:genminorsD} In type $\mathsf{D_n}$:
\begin{enumerate}
  \item For $k =1, \ldots, n-2$, $F_{c^m \omega_k}(t_1, \dots, t_n)$ equals the sum of weights of all collections of vertex-disjoint paths in $\Gamma (\mathsf{D_n}, c)$ with the sources and sinks labeled by~$c^m \cdot [1, k]$;
  \item $F_{c^m \omega_{n-1}}(t_1, \dots, t_n)$ equals the sum of \textbf{square roots} of weights of all collections of \textbf{bundled} vertex-disjoint paths in $\Gamma (\mathsf{D_n}, c)$ with the sources and the sinks labeled by~$c^m \cdot \{1, 2, \ldots, n-1, \barred{n}\}$;
  \item $F_{c^m \omega_{n}}(t_1, \dots, t_n)$ equals the sum of \textbf{square roots} of weights of all collections of \textbf{bundled} vertex-disjoint paths in $\Gamma (\mathsf{D_n}, c)$ with the sources and the sinks labeled by~$c^m \cdot \{1, 2, \ldots, n-1, n\}$.
\end{enumerate}
\end{prop}

The proof will be given in Section~\ref{sec:proofs}, here is an example to illustrate this proposition.

\begin{example}
Type $\mathsf{D_4}$: Let $c = s_1s_2s_3s_4 = (1,2,3,\barred{1},\barred{2},\barred{3}) (4, \barred{4})$, then $c^2 \cdot[1, 2]=\{3, \barred{1}\}$. We have
$$F_{c^2 \omega_2}(t_1, t_2, t_3, t_4) = 1+t_1+t_2+2t_1t_2+t_1t_2t_3+t_1t_2t_4+t_1t_2^2+t_1t_2^2t_3+t_1t_2^2t_4+t_1t_2^2t_3t_4\ \mbox{(see Figure~\ref{fig:D4ex})}.$$

\begin{figure}[h!]
\setlength{\unitlength}{1.2pt}
\begin{center}
\begin{picture}(90,90)(0,-10)
\thicklines
\drawline[90](0,0)(80,0)
\drawline[90](0,10)(80,10)
\drawline[90](0,20)(80,20)
\drawline[90](0,30)(80,30)
\drawline[90](0,40)(80,40)
\drawline[90](0,50)(80,50)
\drawline[90](0,60)(80,60)
\drawline[90](0,70)(80,70)
  \put(0,0){\circle*{2}}
  \put(0,10){\circle*{2}}
  \put(0,20){\circle*{2}}
  \put(0,30){\circle*{2}}
  \put(0,40){\circle*{2}}
  \put(0,50){\circle*{2}}
  \put(0,60){\circle*{2}}
  \put(0,70){\circle*{2}}
  \put(80,0){\circle*{2}}
  \put(80,10){\circle*{2}}
  \put(80,20){\circle*{2}}
  \put(80,30){\circle*{2}}
  \put(80,40){\circle*{2}}
  \put(80,50){\circle*{2}}
  \put(80,60){\circle*{2}}
  \put(80,70){\circle*{2}}
\put(10,0){\circle*{2}}
\put(10,20){\circle*{2}}
\put(10,60){\circle*{2}}
\put(20,10){\circle*{2}}
\put(20,30){\circle*{2}}
\put(20,50){\circle*{2}}
\put(30,20){\circle*{2}}
\put(30,40){\circle*{2}}
\put(30,50){\circle*{2}}
\put(40,20){\circle*{2}}
\put(40,30){\circle*{2}}
\put(50,20){\circle*{2}}
\put(50,40){\circle*{2}}
\put(50,50){\circle*{2}}
\put(60,10){\circle*{2}}
\put(60,30){\circle*{2}}
\put(60,50){\circle*{2}}
\put(70,0){\circle*{2}}
\put(70,20){\circle*{2}}
\put(70,60){\circle*{2}}
  \put(-8,-2){\scriptsize$1$}
  \put(-8,8){\scriptsize$2$}
  \put(-8,18){\scriptsize$3$}
  \put(-8,28){\scriptsize$4$}
  \put(-8,38){\scriptsize$\barred{4}$}
  \put(-8,48){\scriptsize$\barred{3}$}
  \put(-8,58){\scriptsize$\barred{2}$}
  \put(-8,68){\scriptsize$\barred{1}$}
  \put(84,-2){\scriptsize$1$}
  \put(84,8){\scriptsize$2$}
  \put(84,18){\scriptsize$3$}
  \put(84,28){\scriptsize$4$}
  \put(84,38){\scriptsize$\barred{4}$}
  \put(84,48){\scriptsize$\barred{3}$}
  \put(84,58){\scriptsize$\barred{2}$}
  \put(84,68){\scriptsize$\barred{1}$}
\light{
\drawline[90](80,10)(70,0)
\drawline[90](80,70)(50,40)(40,20)
\drawline[90](70,20)(60,10)
\drawline[90](60,30)(50,20)
\drawline[90](50,50)(40,30)
  \put(77,3){\scriptsize$t_1$}
  \put(77,63){\scriptsize$t_1$}
  \put(66,12){\scriptsize$t_2$}
  \put(66,52){\scriptsize$t_2$}
  \put(56,22){\scriptsize$t_3$}
  \put(56,42){\scriptsize$t_3$}
  \put(45,25){\scriptsize$t_4$}
  \put(40,44){\scriptsize$t_4$}
}
\dark{
\drawline[90](40,30)(30,50)
\drawline[90](40,20)(30,40)(0,70)
\drawline[90](30,20)(20,30)
\drawline[90](20,10)(10,20)
\drawline[90](10,0)(0,10)
}
\put(14,-12){$\Gamma (\mathsf{D_4}, s_1s_2s_3s_4)$}
\end{picture}
\begin{picture}(90,90)(0,-10)
\thinlines
\dottedline{2}(0,0)(80,0)
\dottedline{2}(0,10)(80,10)
\dottedline{2}(0,20)(80,20)
\dottedline{2}(0,30)(80,30)
\dottedline{2}(0,40)(80,40)
\dottedline{2}(0,50)(80,50)
\dottedline{2}(0,60)(80,60)
\dottedline{2}(0,70)(80,70)
  \put(0,20){\circle*{2}}
  \put(0,70){\circle*{2}}
  \put(80,20){\circle*{2}}
  \put(80,70){\circle*{2}}
\dottedline{2}(80,70)(50,40)(40,20)(30,40)(0,70)
\dottedline{2}(50,50)(40,30)(30,50)
\dottedline{2}(80,10)(70,0)
\dottedline{2}(70,20)(60,10)
\dottedline{2}(60,30)(50,20)
\dottedline{2}(30,20)(20,30)
\dottedline{2}(20,10)(10,20)
\dottedline{2}(10,0)(0,10)
\thicklines
\drawline[90](0,20)(80,20)
\drawline[90](0,70)(80,70)
\put(18,-12){$\mathrm{wt}: 1$}
\end{picture}
\begin{picture}(90,90)(0,-10)
\thinlines
\dottedline{2}(0,0)(80,0)
\dottedline{2}(0,10)(80,10)
\dottedline{2}(0,20)(80,20)
\dottedline{2}(0,30)(80,30)
\dottedline{2}(0,40)(80,40)
\dottedline{2}(0,50)(80,50)
\dottedline{2}(0,60)(80,60)
\dottedline{2}(0,70)(80,70)
  \put(0,20){\circle*{2}}
  \put(0,70){\circle*{2}}
  \put(80,20){\circle*{2}}
  \put(80,70){\circle*{2}}
\dottedline{2}(80,70)(50,40)(40,20)(30,40)(0,70)
\dottedline{2}(50,50)(40,30)(30,50)
\dottedline{2}(80,10)(70,0)
\dottedline{2}(70,20)(60,10)
\dottedline{2}(60,30)(50,20)
\dottedline{2}(30,20)(20,30)
\dottedline{2}(20,10)(10,20)
\dottedline{2}(10,0)(0,10)
\thicklines
\drawline[90](0,20)(80,20)
\drawline[90](80,70)(70,60)(10,60)(0,70)
\put(18,-12){$\mathrm{wt}: t_1$}
\end{picture}
\begin{picture}(90,90)(0,-10)
\thinlines
\dottedline{2}(0,0)(80,0)
\dottedline{2}(0,10)(80,10)
\dottedline{2}(0,20)(80,20)
\dottedline{2}(0,30)(80,30)
\dottedline{2}(0,40)(80,40)
\dottedline{2}(0,50)(80,50)
\dottedline{2}(0,60)(80,60)
\dottedline{2}(0,70)(80,70)
  \put(0,20){\circle*{2}}
  \put(0,70){\circle*{2}}
  \put(80,20){\circle*{2}}
  \put(80,70){\circle*{2}}
\dottedline{2}(80,70)(50,40)(40,20)(30,40)(0,70)
\dottedline{2}(50,50)(40,30)(30,50)
\dottedline{2}(80,10)(70,0)
\dottedline{2}(70,20)(60,10)
\dottedline{2}(60,30)(50,20)
\dottedline{2}(30,20)(20,30)
\dottedline{2}(20,10)(10,20)
\dottedline{2}(10,0)(0,10)
\thicklines
\drawline[90](0,20)(80,20)
\drawline[90](80,70)(60,50)(20,50)(0,70)
\put(18,-12){$\mathrm{wt}: t_1t_2$}
\end{picture}
\begin{picture}(90,100)(0,-10)
\thinlines
\dottedline{2}(0,0)(80,0)
\dottedline{2}(0,10)(80,10)
\dottedline{2}(0,20)(80,20)
\dottedline{2}(0,30)(80,30)
\dottedline{2}(0,40)(80,40)
\dottedline{2}(0,50)(80,50)
\dottedline{2}(0,60)(80,60)
\dottedline{2}(0,70)(80,70)
  \put(0,20){\circle*{2}}
  \put(0,70){\circle*{2}}
  \put(80,20){\circle*{2}}
  \put(80,70){\circle*{2}}
\dottedline{2}(80,70)(50,40)(40,20)(30,40)(0,70)
\dottedline{2}(50,50)(40,30)(30,50)
\dottedline{2}(80,10)(70,0)
\dottedline{2}(70,20)(60,10)
\dottedline{2}(60,30)(50,20)
\dottedline{2}(30,20)(20,30)
\dottedline{2}(20,10)(10,20)
\dottedline{2}(10,0)(0,10)
\thicklines
\drawline[90](0,20)(80,20)
\drawline[90](80,70)(50,40)(30,40)(0,70)
\put(18,-12){$\mathrm{wt}: t_1t_2t_3$}
\end{picture}
\begin{picture}(90,100)(0,-10)
\thinlines
\dottedline{2}(0,0)(80,0)
\dottedline{2}(0,10)(80,10)
\dottedline{2}(0,20)(80,20)
\dottedline{2}(0,30)(80,30)
\dottedline{2}(0,40)(80,40)
\dottedline{2}(0,50)(80,50)
\dottedline{2}(0,60)(80,60)
\dottedline{2}(0,70)(80,70)
  \put(0,20){\circle*{2}}
  \put(0,70){\circle*{2}}
  \put(80,20){\circle*{2}}
  \put(80,70){\circle*{2}}
\dottedline{2}(80,70)(50,40)(40,20)(30,40)(0,70)
\dottedline{2}(50,50)(40,30)(30,50)
\dottedline{2}(80,10)(70,0)
\dottedline{2}(70,20)(60,10)
\dottedline{2}(60,30)(50,20)
\dottedline{2}(30,20)(20,30)
\dottedline{2}(20,10)(10,20)
\dottedline{2}(10,0)(0,10)
\thicklines
\drawline[90](0,20)(80,20)
\drawline[90](80,70)(60,50)(50,50)(40,30)(30,50)(20,50)(0,70)
\put(18,-12){$\mathrm{wt}: t_1t_2t_4$}
\end{picture}
\begin{picture}(90,100)(0,-10)
\thinlines
\dottedline{2}(0,0)(80,0)
\dottedline{2}(0,10)(80,10)
\dottedline{2}(0,20)(80,20)
\dottedline{2}(0,30)(80,30)
\dottedline{2}(0,40)(80,40)
\dottedline{2}(0,50)(80,50)
\dottedline{2}(0,60)(80,60)
\dottedline{2}(0,70)(80,70)
  \put(0,20){\circle*{2}}
  \put(0,70){\circle*{2}}
  \put(80,20){\circle*{2}}
  \put(80,70){\circle*{2}}
\dottedline{2}(80,70)(50,40)(40,20)(30,40)(0,70)
\dottedline{2}(50,50)(40,30)(30,50)
\dottedline{2}(80,10)(70,0)
\dottedline{2}(70,20)(60,10)
\dottedline{2}(60,30)(50,20)
\dottedline{2}(30,20)(20,30)
\dottedline{2}(20,10)(10,20)
\dottedline{2}(10,0)(0,10)
\thicklines
\drawline[90](80,20)(70,20)(60,10)(20,10)(10,20)(0,20)
\drawline[90](0,70)(80,70)
\put(18,-12){$\mathrm{wt}: t_2$}
\end{picture}
\begin{picture}(90,100)(0,-10)
\thinlines
\dottedline{2}(0,0)(80,0)
\dottedline{2}(0,10)(80,10)
\dottedline{2}(0,20)(80,20)
\dottedline{2}(0,30)(80,30)
\dottedline{2}(0,40)(80,40)
\dottedline{2}(0,50)(80,50)
\dottedline{2}(0,60)(80,60)
\dottedline{2}(0,70)(80,70)
  \put(0,20){\circle*{2}}
  \put(0,70){\circle*{2}}
  \put(80,20){\circle*{2}}
  \put(80,70){\circle*{2}}
\dottedline{2}(80,70)(50,40)(40,20)(30,40)(0,70)
\dottedline{2}(50,50)(40,30)(30,50)
\dottedline{2}(80,10)(70,0)
\dottedline{2}(70,20)(60,10)
\dottedline{2}(60,30)(50,20)
\dottedline{2}(30,20)(20,30)
\dottedline{2}(20,10)(10,20)
\dottedline{2}(10,0)(0,10)
\thicklines
\drawline[90](80,20)(70,20)(60,10)(20,10)(10,20)(0,20)
\drawline[90](80,70)(70,60)(10,60)(0,70)
\put(18,-12){$\mathrm{wt}: t_1t_2$}
\end{picture}
\begin{picture}(90,100)(0,-10)
\thinlines
\dottedline{2}(0,0)(80,0)
\dottedline{2}(0,10)(80,10)
\dottedline{2}(0,20)(80,20)
\dottedline{2}(0,30)(80,30)
\dottedline{2}(0,40)(80,40)
\dottedline{2}(0,50)(80,50)
\dottedline{2}(0,60)(80,60)
\dottedline{2}(0,70)(80,70)
  \put(0,20){\circle*{2}}
  \put(0,70){\circle*{2}}
  \put(80,20){\circle*{2}}
  \put(80,70){\circle*{2}}
\dottedline{2}(80,70)(50,40)(40,20)(30,40)(0,70)
\dottedline{2}(50,50)(40,30)(30,50)
\dottedline{2}(80,10)(70,0)
\dottedline{2}(70,20)(60,10)
\dottedline{2}(60,30)(50,20)
\dottedline{2}(30,20)(20,30)
\dottedline{2}(20,10)(10,20)
\dottedline{2}(10,0)(0,10)
\thicklines
\drawline[90](80,20)(70,20)(60,10)(20,10)(10,20)(0,20)
\drawline[90](80,70)(60,50)(20,50)(0,70)
\put(18,-12){$\mathrm{wt}: t_1t_2^2$}
\end{picture}
\begin{picture}(90,100)(0,-10)
\thinlines
\dottedline{2}(0,0)(80,0)
\dottedline{2}(0,10)(80,10)
\dottedline{2}(0,20)(80,20)
\dottedline{2}(0,30)(80,30)
\dottedline{2}(0,40)(80,40)
\dottedline{2}(0,50)(80,50)
\dottedline{2}(0,60)(80,60)
\dottedline{2}(0,70)(80,70)
  \put(0,20){\circle*{2}}
  \put(0,70){\circle*{2}}
  \put(80,20){\circle*{2}}
  \put(80,70){\circle*{2}}
\dottedline{2}(80,70)(50,40)(40,20)(30,40)(0,70)
\dottedline{2}(50,50)(40,30)(30,50)
\dottedline{2}(80,10)(70,0)
\dottedline{2}(70,20)(60,10)
\dottedline{2}(60,30)(50,20)
\dottedline{2}(30,20)(20,30)
\dottedline{2}(20,10)(10,20)
\dottedline{2}(10,0)(0,10)
\thicklines
\drawline[90](80,20)(70,20)(60,10)(20,10)(10,20)(0,20)
\drawline[90](80,70)(50,40)(30,40)(0,70)
\put(18,-12){$\mathrm{wt}: t_1t_2^2t_3$}
\end{picture}
\begin{picture}(90,100)(0,-10)
\thinlines
\dottedline{2}(0,0)(80,0)
\dottedline{2}(0,10)(80,10)
\dottedline{2}(0,20)(80,20)
\dottedline{2}(0,30)(80,30)
\dottedline{2}(0,40)(80,40)
\dottedline{2}(0,50)(80,50)
\dottedline{2}(0,60)(80,60)
\dottedline{2}(0,70)(80,70)
  \put(0,20){\circle*{2}}
  \put(0,70){\circle*{2}}
  \put(80,20){\circle*{2}}
  \put(80,70){\circle*{2}}
\dottedline{2}(80,70)(50,40)(40,20)(30,40)(0,70)
\dottedline{2}(50,50)(40,30)(30,50)
\dottedline{2}(80,10)(70,0)
\dottedline{2}(70,20)(60,10)
\dottedline{2}(60,30)(50,20)
\dottedline{2}(30,20)(20,30)
\dottedline{2}(20,10)(10,20)
\dottedline{2}(10,0)(0,10)
\thicklines
\drawline[90](80,20)(70,20)(60,10)(20,10)(10,20)(0,20)
\drawline[90](80,70)(60,50)(50,50)(40,30)(30,50)(20,50)(0,70)
\put(18,-12){$\mathrm{wt}: t_1t_2^2t_4$}
\end{picture}
\begin{picture}(90,100)(0,-10)
\thinlines
\dottedline{2}(0,0)(80,0)
\dottedline{2}(0,10)(80,10)
\dottedline{2}(0,20)(80,20)
\dottedline{2}(0,30)(80,30)
\dottedline{2}(0,40)(80,40)
\dottedline{2}(0,50)(80,50)
\dottedline{2}(0,60)(80,60)
\dottedline{2}(0,70)(80,70)
  \put(0,20){\circle*{2}}
  \put(0,70){\circle*{2}}
  \put(80,20){\circle*{2}}
  \put(80,70){\circle*{2}}
\dottedline{2}(80,70)(50,40)(40,20)(30,40)(0,70)
\dottedline{2}(50,50)(40,30)(30,50)
\dottedline{2}(80,10)(70,0)
\dottedline{2}(70,20)(60,10)
\dottedline{2}(60,30)(50,20)
\dottedline{2}(30,20)(20,30)
\dottedline{2}(20,10)(10,20)
\dottedline{2}(10,0)(0,10)
\thicklines
\drawline[90](80,20)(70,20)(60,10)(20,10)(10,20)(0,20)
\drawline[90](80,70)(50,40)(40,20)(30,40)(0,70)
\put(18,-12){$\mathrm{wt}: t_1t_2^2t_3t_4$}
\end{picture}
\end{center}
\caption{$\Gamma (D_4, s_1s_2s_3s_4)$, $F_{c^2 \omega_2} (t_1,t_2,t_3,t_4)$.} \label{fig:D4ex}
\end{figure}
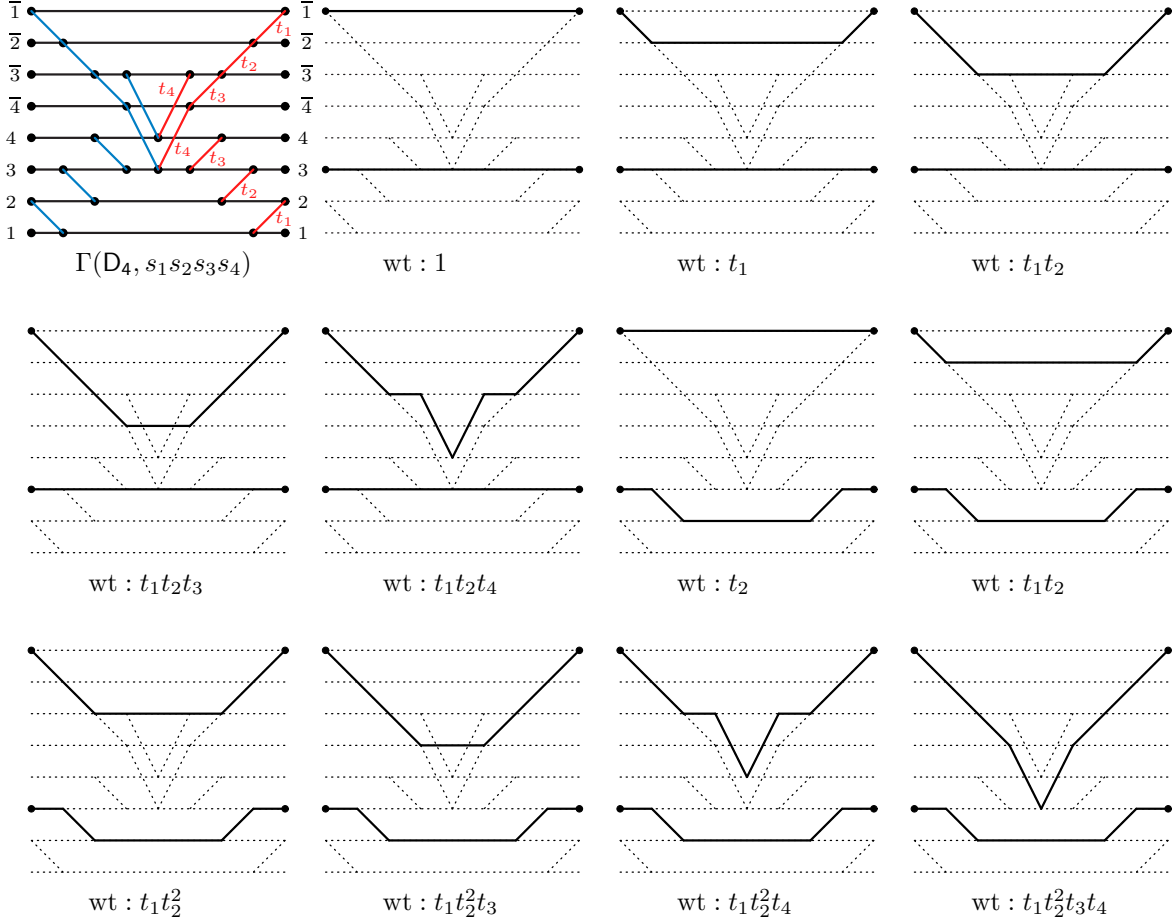

We also have  $c^2 \cdot \{1, 2, 3, \barred{4}\}=\{3, \barred{4}, \barred{2}, \barred{1}\}$, hence $F_{c^2 \omega_3}(t_1, t_2, t_3, t_4) = 1 + t_2 + t_2t_4$ (see Figure~\ref{fig:D4spinex}). Remember that in this case we require bundled families of paths and only square roots of their weights contribute to the $F$-polynomial.

\begin{figure}[h!]
\setlength{\unitlength}{1.2pt}
\begin{center}
\begin{picture}(90,80)(0,-10)
\thinlines
\dottedline{2}(0,0)(80,0)
\dottedline{2}(0,10)(80,10)
\dottedline{2}(0,20)(80,20)
\dottedline{2}(0,30)(80,30)
\dottedline{2}(0,40)(80,40)
\dottedline{2}(0,50)(80,50)
\dottedline{2}(0,60)(80,60)
\dottedline{2}(0,70)(80,70)
  \put(0,20){\circle*{2}}
  \put(0,40){\circle*{2}}
  \put(0,60){\circle*{2}}
  \put(0,70){\circle*{2}}
  \put(80,20){\circle*{2}}
  \put(80,40){\circle*{2}}
  \put(80,60){\circle*{2}}
  \put(80,70){\circle*{2}}
\dottedline{2}(80,70)(50,40)(40,20)(30,40)(0,70)
\dottedline{2}(50,50)(40,30)(30,50)
\dottedline{2}(80,10)(70,0)
\dottedline{2}(70,20)(60,10)
\dottedline{2}(60,30)(50,20)
\dottedline{2}(30,20)(20,30)
\dottedline{2}(20,10)(10,20)
\dottedline{2}(10,0)(0,10)
\thicklines
\drawline[90](80,20)(0,20)
\drawline[90](80,40)(0,40)
\drawline[90](80,60)(0,60)
\drawline[90](80,70)(0,70)
\put(18,-12){$\mathrm{wt}: 1$}
\end{picture}
\begin{picture}(90,80)(0,-10)
\thinlines
\dottedline{2}(0,0)(80,0)
\dottedline{2}(0,10)(80,10)
\dottedline{2}(0,20)(80,20)
\dottedline{2}(0,30)(80,30)
\dottedline{2}(0,40)(80,40)
\dottedline{2}(0,50)(80,50)
\dottedline{2}(0,60)(80,60)
\dottedline{2}(0,70)(80,70)
  \put(0,20){\circle*{2}}
  \put(0,40){\circle*{2}}
  \put(0,60){\circle*{2}}
  \put(0,70){\circle*{2}}
  \put(80,20){\circle*{2}}
  \put(80,40){\circle*{2}}
  \put(80,60){\circle*{2}}
  \put(80,70){\circle*{2}}
\dottedline{2}(80,70)(50,40)(40,20)(30,40)(0,70)
\dottedline{2}(50,50)(40,30)(30,50)
\dottedline{2}(80,10)(70,0)
\dottedline{2}(70,20)(60,10)
\dottedline{2}(60,30)(50,20)
\dottedline{2}(30,20)(20,30)
\dottedline{2}(20,10)(10,20)
\dottedline{2}(10,0)(0,10)
\thicklines
\drawline[90](80,20)(70,20)(60,10)(20,10)(10,20)(0,20)
\drawline[90](80,40)(0,40)
\drawline[90](80,60)(70,60)(60,50)(20,50)(10,60)(0,60)
\drawline[90](80,70)(0,70)
\put(18,-12){$\mathrm{wt}: t_2^2$}
\end{picture}
\begin{picture}(90,80)(0,-10)
\thinlines
\dottedline{2}(0,0)(80,0)
\dottedline{2}(0,10)(80,10)
\dottedline{2}(0,20)(80,20)
\dottedline{2}(0,30)(80,30)
\dottedline{2}(0,40)(80,40)
\dottedline{2}(0,50)(80,50)
\dottedline{2}(0,60)(80,60)
\dottedline{2}(0,70)(80,70)
  \put(0,20){\circle*{2}}
  \put(0,40){\circle*{2}}
  \put(0,60){\circle*{2}}
  \put(0,70){\circle*{2}}
  \put(80,20){\circle*{2}}
  \put(80,40){\circle*{2}}
  \put(80,60){\circle*{2}}
  \put(80,70){\circle*{2}}
\dottedline{2}(80,70)(50,40)(40,20)(30,40)(0,70)
\dottedline{2}(50,50)(40,30)(30,50)
\dottedline{2}(80,10)(70,0)
\dottedline{2}(70,20)(60,10)
\dottedline{2}(60,30)(50,20)
\dottedline{2}(30,20)(20,30)
\dottedline{2}(20,10)(10,20)
\dottedline{2}(10,0)(0,10)
\thicklines
\drawline[90](80,20)(70,20)(60,10)(20,10)(10,20)(0,20)
\drawline[90](80,40)(50,40)(40,20)(30,40)(0,40)
\drawline[90](80,60)(70,60)(60,50)(50,50)(40,30)(30,50)(20,50)(10,60)(0,60)
\drawline[90](80,70)(0,70)
\put(18,-12){$\mathrm{wt}: t_2^2t_4^2$}
\end{picture}
\end{center}
\caption{$\Gamma (D_4, s_1s_2s_3s_4)$, $F_{c^2 \omega_3} (t_1,t_2,t_3,t_4)$.} \label{fig:D4spinex}
\end{figure}
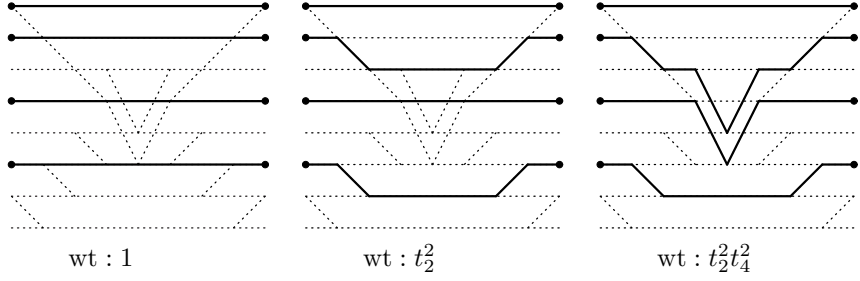

\end{example}


\textbf{Type $\mathsf{B_n}$} ($n \geq 2$) :
For each $i \in [1, n] \cup [\barred{1}, \barred{n}]$, the elementary chip corresponding to $x_{i}(t)$ is shown in Figure~\ref{fig:bnchipsnsp}.
\begin{figure}[ht]
\setlength{\unitlength}{1.2pt}
\begin{center}
\begin{picture}(63,120)(0,-15)
\thicklines
\put(18,0){\line(1,0){20}} \put(18,20){\line(1,0){20}} \put(18,40){\line(1,0){20}} \put(18,60){\line(1,0){20}}
\put(18,80){\line(1,0){20}} \put(18,100){\line(1,0){20}}
  \put(18,0){\circle*{2.5}}
  \put(18,20){\circle*{2.5}}
  \put(18,40){\circle*{2.5}}
  \put(18,60){\circle*{2.5}}
  \put(18,80){\circle*{2.5}}
  \put(18,100){\circle*{2.5}}
  \put(38,0){\circle*{2.5}}
  \put(38,20){\circle*{2.5}}
  \put(38,40){\circle*{2.5}}
  \put(38,60){\circle*{2.5}}
  \put(38,80){\circle*{2.5}}
  \put(38,100){\circle*{2.5}}
  \put(4,-2){\scriptsize${1}$}
  \put(4,18){\scriptsize${i}$}
  \put(0,38){\scriptsize${i+1}$}
  \put(0,58){\scriptsize${\overline{i+1}}$}
  \put(4,78){\scriptsize${\overline{i}}$}
  \put(4,98){\scriptsize${\overline{1}}$}
  \dottedline{3}(6,6)(6,13)
  \dottedline{3}(6,46)(6,53)
  \dottedline{3}(6,86)(6,93)
  \dottedline{3}(28,5)(28,15)
  \dottedline{3}(28,45)(28,55)
  \dottedline{3}(28,85)(28,95)
\dark{ \put(18,40){\line(1,-1){20}} \put(18,80){\line(1,-1){20}} }
  \put(22,26){\footnotesize$t$}
  \put(22,66){\footnotesize$t$}
  \put(17,-17){$x_{\barred{i}} (t)$}
\end{picture}
\begin{picture}(50,120)(0,-15)
\thicklines
\put(5,0){\line(1,0){20}} \put(5,20){\line(1,0){20}} \put(5,40){\line(1,0){20}} \put(5,60){\line(1,0){20}}
\put(5,80){\line(1,0){20}} \put(5,100){\line(1,0){20}}
  \put(5,0){\circle*{2.5}}
  \put(5,20){\circle*{2.5}}
  \put(5,40){\circle*{2.5}}
  \put(5,60){\circle*{2.5}}
  \put(5,80){\circle*{2.5}}
  \put(5,100){\circle*{2.5}}
  \put(25,0){\circle*{2.5}}
  \put(25,20){\circle*{2.5}}
  \put(25,40){\circle*{2.5}}
  \put(25,60){\circle*{2.5}}
  \put(25,80){\circle*{2.5}}
  \put(25,100){\circle*{2.5}}
  \dottedline{3}(15,5)(15,15)
  \dottedline{3}(15,45)(15,55)
  \dottedline{3}(15,85)(15,95)
\light{ \put(5,20){\line(1,1){20}} \put(5,60){\line(1,1){20}} }
  \put(18,26){\footnotesize$t$}
  \put(18,66){\footnotesize$t$}
  \put(4,-17){$x_i (t)$}
\end{picture}
\begin{picture}(50,120)(0,-15)
\thicklines
\put(5,0){\line(1,0){20}} \put(5,30){\line(1,0){20}} \put(5,50){\line(1,0){20}} \put(5,70){\line(1,0){20}} \put(5,100){\line(1,0){20}}
  \put(5,0){\circle*{2.5}}
  \put(5,30){\circle*{2.5}}
  \put(5,50){\circle*{2.5}}
  \put(5,70){\circle*{2.5}}
  \put(5,100){\circle*{2.5}}
  \put(25,0){\circle*{2.5}}
  \put(25,30){\circle*{2.5}}
  \put(25,50){\circle*{2.5}}
  \put(25,70){\circle*{2.5}}
  \put(25,100){\circle*{2.5}}
  \dottedline{3}(15,5)(15,25)
  \dottedline{3}(15,75)(15,95)
\dark{ \put(5,70){\line(1,-1){20}} \put(5,70){\line(1,-2){20}} \put(5,50){\line(1,-1){20}} }
  \put(17,59){\footnotesize$\sqrt{2}t$}
  \put(0,35){\footnotesize$\sqrt{2}t$}  \put(20,41){\footnotesize$t^2$}
  \put(4,-17){$x_{\barred{n}} (t)$}
\end{picture}
\begin{picture}(50,120)(0,-15)
\thicklines
\put(5,0){\line(1,0){20}} \put(5,30){\line(1,0){20}} \put(5,50){\line(1,0){20}} \put(5,70){\line(1,0){20}} \put(5,100){\line(1,0){20}}
  \put(5,0){\circle*{2.5}}
  \put(5,30){\circle*{2.5}}
  \put(5,50){\circle*{2.5}}
  \put(5,70){\circle*{2.5}}
  \put(5,100){\circle*{2.5}}
  \put(25,0){\circle*{2.5}}
  \put(25,30){\circle*{2.5}}
  \put(25,50){\circle*{2.5}}
  \put(25,70){\circle*{2.5}}
  \put(25,100){\circle*{2.5}}
  \dottedline{3}(15,5)(15,25)
  \dottedline{3}(15,75)(15,95)
  \put(35,-2){\scriptsize${1}$}
  \put(35,28){\scriptsize${n}$}
  \put(35,48){\scriptsize${0}$}
  \put(35,68){\scriptsize${\overline{n}}$}
  \put(35,98){\scriptsize${\overline{1}}$}
  \dottedline{3}(36,8)(36,25)
  \dottedline{3}(36,78)(36,95)
\light{ \put(5,30){\line(1,1){20}} \put(5,30){\line(1,2){20}} \put(5,50){\line(1,1){20}} }
  \put(0,60){\footnotesize$\sqrt{2}t$} \put(21,55){\footnotesize$t^2$}
  \put(17,37){\footnotesize$\sqrt{2}t$}
  \put(4,-17){$x_n (t)$}
\end{picture}
\end{center}
\caption{Elementary chips in $\Gamma (\mathsf{B_n}, \ii)$ ($i = 1, \ldots, n-1$).} \label{fig:bnchipsnsp}
\end{figure}
In each chip, the vertices consist of all the endpoints of the $2n+1$ horizontal edges. All the edges are oriented from right to left, we number the horizontal levels from bottom to top in the order $1, \ldots, n, 0, \barred{n}, \ldots, \barred{1}$. In the chips corresponding to $x_{\barred{n}}(t)$ and $x_n(t)$, the intersections of the diagonal edges and the horizontal edges in the middle of the diagonal edges on the horizontal level $0$ are \emph{not} vertices.
The numbering of the horizontal levels for the first (resp., last) two chips in Figure~\ref{fig:bnchipsnsp} is shown on the left (resp., right) of the figure. All unlabeled edges have weight 1.
The directed graph $\Gamma (\mathsf{B_n}, c)$ associated with $c = s_{i_1} \cdots s_{i_n}$ is constructed as a concatenation of elementary chips $x_{\barred{i_i}}(1), \ldots, x_{\barred{i_n}}(1), x_{i_n}(t_{i_n}), \ldots, x_{i_1}(t_{i_1})$ (in this order). We number the $2n+1$ sources and the $2n+1$ sinks of the graph $\Gamma (\mathsf{B_n}, c)$ bottom-to-top in the order $1, \ldots, n, 0, \barred{n}, \ldots, \barred{1}$.

To finish the type $\mathsf{B_n}$ case, we need to introduce another graph $\Gamma_{\mathrm{S}} (\mathsf{B_n}, c)$ (it corresponds to the spin representation). For each $i \in [1, n] \cup [\barred{1}, \barred{n}]$, the elementary chip corresponding to $x_{i}(t)$ in $\Gamma_{\mathrm{S}} (\mathsf{B_n}, c)$ is shown in Figure~\ref{fig:bnchipssp}. The vertices for each elementary chip consist of all the endpoints of the $2n$ horizontal edges. We label the $2n$ horizontal levels from bottom to top by $1, \ldots, n, \barred{n}, \ldots, \barred{1}$. All the edges are oriented from right to left with their weights shown in the figure, all unlabeled edges have weight 1.
\begin{figure}[ht]
\setlength{\unitlength}{1.2pt}
\begin{center}
\begin{picture}(63,120)(0,-15)
\thicklines
\put(18,0){\line(1,0){20}} \put(18,20){\line(1,0){20}} \put(18,40){\line(1,0){20}} \put(18,60){\line(1,0){20}}
\put(18,80){\line(1,0){20}} \put(18,100){\line(1,0){20}}
  \put(18,0){\circle*{2.5}}
  \put(18,20){\circle*{2.5}}
  \put(18,40){\circle*{2.5}}
  \put(18,60){\circle*{2.5}}
  \put(18,80){\circle*{2.5}}
  \put(18,100){\circle*{2.5}}
  \put(38,0){\circle*{2.5}}
  \put(38,20){\circle*{2.5}}
  \put(38,40){\circle*{2.5}}
  \put(38,60){\circle*{2.5}}
  \put(38,80){\circle*{2.5}}
  \put(38,100){\circle*{2.5}}
  \put(4,-2){\scriptsize${1}$}
  \put(4,18){\scriptsize${i}$}
  \put(0,38){\scriptsize${i+1}$}
  \put(0,58){\scriptsize${\overline{i+1}}$}
  \put(4,78){\scriptsize${\overline{i}}$}
  \put(4,98){\scriptsize${\overline{1}}$}
  \dottedline{3}(6,6)(6,13)
  \dottedline{3}(6,46)(6,53)
  \dottedline{3}(6,86)(6,93)
  \dottedline{3}(28,5)(28,15)
  \dottedline{3}(28,45)(28,55)
  \dottedline{3}(28,85)(28,95)
\dark{ \put(18,40){\line(1,-1){20}} \put(18,80){\line(1,-1){20}} }
  \put(20,26){\footnotesize$t$}
  \put(20,66){\footnotesize$t$}
  \put(17,-17){$x_{\barred{i}} (t)$}
\end{picture}
\begin{picture}(50,120)(0,-15)
\thicklines
\put(5,0){\line(1,0){20}} \put(5,20){\line(1,0){20}} \put(5,40){\line(1,0){20}} \put(5,60){\line(1,0){20}}
\put(5,80){\line(1,0){20}} \put(5,100){\line(1,0){20}}
  \put(5,0){\circle*{2.5}}
  \put(5,20){\circle*{2.5}}
  \put(5,40){\circle*{2.5}}
  \put(5,60){\circle*{2.5}}
  \put(5,80){\circle*{2.5}}
  \put(5,100){\circle*{2.5}}
  \put(25,0){\circle*{2.5}}
  \put(25,20){\circle*{2.5}}
  \put(25,40){\circle*{2.5}}
  \put(25,60){\circle*{2.5}}
  \put(25,80){\circle*{2.5}}
  \put(25,100){\circle*{2.5}}
  \dottedline{3}(15,5)(15,15)
  \dottedline{3}(15,45)(15,55)
  \dottedline{3}(15,85)(15,95)
\light{ \put(5,20){\line(1,1){20}} \put(5,60){\line(1,1){20}} }
  \put(18,26){\footnotesize$t$}
  \put(18,66){\footnotesize$t$}
  \put(4,-17){$x_i (t)$}
\end{picture}
\begin{picture}(50,120)(0,-15)
\thicklines
\put(5,0){\line(1,0){20}}  \put(5,40){\line(1,0){20}} \put(5,60){\line(1,0){20}} \put(5,100){\line(1,0){20}}
  \put(5,0){\circle*{2.5}}
  \put(5,40){\circle*{2.5}}
  \put(5,60){\circle*{2.5}}
  \put(5,100){\circle*{2.5}}
  \put(25,0){\circle*{2.5}}
  \put(25,40){\circle*{2.5}}
  \put(25,60){\circle*{2.5}}
  \put(25,100){\circle*{2.5}}
  \dottedline{3}(15,5)(15,35)
  \dottedline{3}(15,65)(15,95)
\dark{ \put(5,60){\line(1,-1){20}} }
  \put(8,44){\footnotesize$t^2$}
  \put(4,-17){$x_{\barred{n}} (t)$}
\end{picture}
\begin{picture}(50,120)(0,-15)
\thicklines
\put(5,0){\line(1,0){20}}  \put(5,40){\line(1,0){20}} \put(5,60){\line(1,0){20}}
 \put(5,100){\line(1,0){20}}
  \put(5,0){\circle*{2.5}}
  \put(5,40){\circle*{2.5}}
  \put(5,60){\circle*{2.5}}
  \put(5,100){\circle*{2.5}}
  \put(25,0){\circle*{2.5}}
  \put(25,40){\circle*{2.5}}
  \put(25,60){\circle*{2.5}}
  \put(25,100){\circle*{2.5}}
  \put(35,-2){\scriptsize${1}$}
  \put(35,38){\scriptsize${n}$}
  \put(35,58){\scriptsize${\overline{n}}$}
  \put(35,98){\scriptsize${\overline{1}}$}
  \dottedline{3}(15,5)(15,35)
  \dottedline{3}(15,65)(15,95)
  \dottedline{3}(36,8)(36,35)
  \dottedline{3}(36,68)(36,95)
\light{ \put(5,40){\line(1,1){20}} }
  \put(17,45){\footnotesize$t^2$}
  \put(4,-17){$x_n (t)$}
\end{picture}
\end{center}
\caption{Elementary chips in $\Gamma_{\mathrm{S}} (\mathsf{B_n}, c)$ ($i = 1, \ldots, n-1$).} \label{fig:bnchipssp}
\end{figure}

The directed graph $\Gamma_{\mathrm{S}} (\mathsf{B_n}, c)$ associated with $c = s_{i_1} \cdots s_{i_n}$ is constructed as a concatenation of elementary chips $x_{\barred{i_i}}(1), \ldots, x_{\barred{i_n}}(1), x_{i_n}(t_{i_n}), \ldots, x_{i_1}(t_{i_1})$ (in this order). We number the $2n$ sources and the $2n$ sinks of the graph $\Gamma_{\mathrm{S}} (\mathsf{B_n}, c)$ bottom-to-top in the order $1, \ldots, n, \barred{n}, \ldots, \barred{1}$.

As before, we call a family of paths in $\Gamma_{\mathrm{S}} (\mathsf{B_n}, c)$ \emph{bundled} if within each elementary chip that corresponds to $x_{i}(t)$, for $i \in [1, n-1] \cup [\barred{1}, \barred{n-1}]$, either both of the diagonal edges belong to the family of paths, or neither belong to the family of paths. We will only need the bundled families of vertex-disjoint paths in $\Gamma_{\mathrm{S}} (\mathsf{B_n}, c)$.

The Weyl group of type $\mathsf{B_n}$ acts on the index set $[1, n] \cup [\barred{1}, \barred{n}]$ by permutations and ``bar'' changes. When written as permutations on $[1, n] \cup [\barred{1}, \barred{n}]$, the simple reflections are $s_i = (i, i+1)(\barred{i+1}, \barred{i})$ for $i = 1, \ldots, n-1$, and $s_n = (n, \overline{n})$. The definition of the weight of paths is the same as before. Then the $F$-polynomials in type $\mathsf{B_n}$ are computed as follows:

\begin{prop}
\label{pr:genminorsB}\mbox{ }
\begin{enumerate}
  \item For $k =1, \ldots, n-1$, $F_{c^m \omega_k}(t_1, \dots, t_n)$ equals the sum of weights of all collections of vertex-disjoint paths in $\Gamma (\mathsf{B_n}, c)$ with the sources and sinks labeled by~$c^m \cdot [1, k]$;
  \item $F_{c^m \omega_n}(t_1, \dots, t_n)$ equals the sum of \textbf{square roots} of weights of all collections of \textbf{bundled} vertex-disjoint paths in $\Gamma_{\mathrm{S}} (\mathsf{B_n}, c)$ with the sources and sinks labeled by~$c^m \cdot [1, n]$.
\end{enumerate}
\end{prop}

\begin{example}
Type $\mathsf{B_2}$: Let $c = s_2s_1 = (2, 1, \barred{2}, \barred{1})$, then $c \cdot[1]=\{\barred{2}\}$ and $c^2\cdot[1,2]=\{\barred{2},\barred{1}\}$. We have
$$F_{c \omega_1}(t_1, t_2) = 1 + 2t_2 + t_2^2 + t_1 t_2^2 \ \ \  \mbox{ and } \ \ \  F_{c^2 \omega_2}(t_1, t_2) = 1 + t_2 + t_1t_2\,.$$

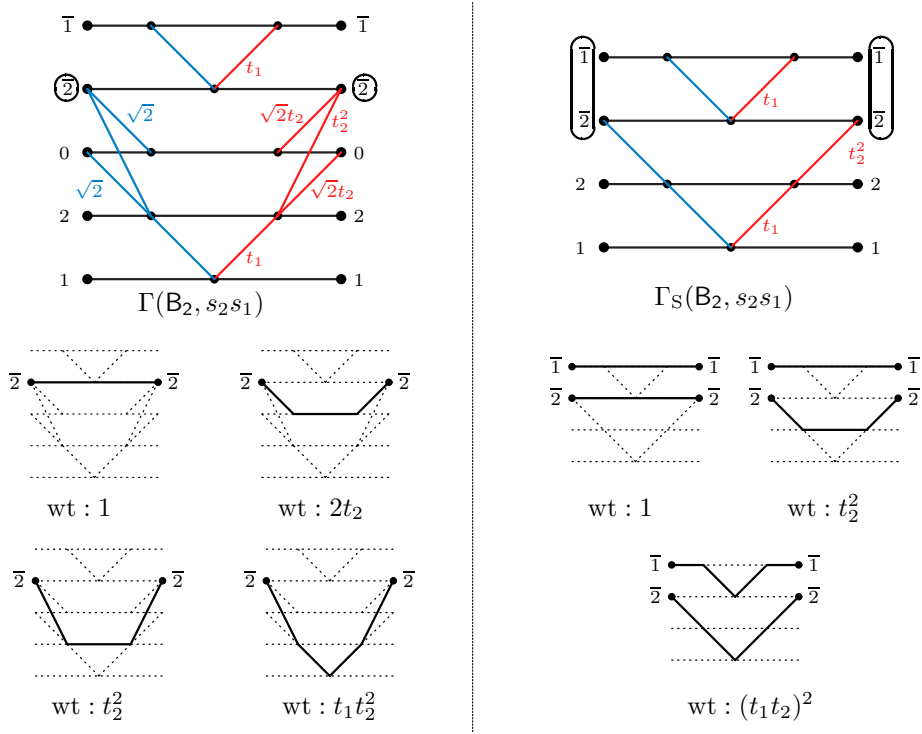
\begin{figure}[h!]
\setlength{\unitlength}{1.2pt}
\begin{center}
\begin{picture}(140,90)(-20,-5)
\thicklines
\drawline[90](0,0)(80,0)
\drawline[90](0,20)(80,20)
\drawline[90](0,40)(80,40)
\drawline[90](0,60)(80,60)
\drawline[90](0,80)(80,80)
  \put(0,0){\circle*{2.5}}
  \put(0,20){\circle*{2.5}}
  \put(0,40){\circle*{2.5}}
  \put(0,60){\circle*{2.5}}
  \put(0,80){\circle*{2.5}}
  \put(80,0){\circle*{2.5}}
  \put(80,20){\circle*{2.5}}
  \put(80,40){\circle*{2.5}}
  \put(80,60){\circle*{2.5}}
  \put(80,80){\circle*{2.5}}
\put(20,20){\circle*{2}}
\put(20,40){\circle*{2}}
\put(20,80){\circle*{2}}
\put(40,0){\circle*{2}}
\put(40,60){\circle*{2}}
\put(60,20){\circle*{2}}
\put(60,40){\circle*{2}}
\put(60,80){\circle*{2}}
  \put(-9,-2){\scriptsize${1}$}
  \put(-9,18){\scriptsize${2}$}
  \put(-9,38){\scriptsize${0}$}
  \put(-9,58){\scriptsize${\barred{2}}$}
  \put(-9,78){\scriptsize${\barred{1}}$}
  \put(84,-2){\scriptsize${1}$}
  \put(84,18){\scriptsize${2}$}
  \put(84,38){\scriptsize${0}$}
  \put(84,58){\scriptsize${\barred{2}}$}
  \put(84,78){\scriptsize${\barred{1}}$}
  \put(87,60){\oval(7,8)}
  \put(-7,60){\oval(7,8)}
\thicklines
\light{
\drawline[90](80,60)(60,40)
\drawline[90](80,60)(60,20)
\drawline[90](80,40)(60,20)
\drawline[90](60,80)(40,60)
\drawline[90](60,20)(40,0)
  \put(50,65){\scriptsize$t_1$}
  \put(50,5){\scriptsize$t_1$}
  \put(54,50){\scriptsize$\sqrt{2}t_2$}
  \put(77,48){\scriptsize$t_2^2$}
  \put(70,26){\scriptsize$\sqrt{2}t_2$}
}
\dark{
\drawline[90](40,60)(20,80)
\drawline[90](40,0)(20,20)
\drawline[90](20,20)(0,40)
\drawline[90](20,20)(0,60)
\drawline[90](20,40)(0,60)
  \put(12,50){\scriptsize$\sqrt{2}$}
  \put(-4,26){\scriptsize$\sqrt{2}$}
}
\put(16,-10){$\Gamma (\mathsf{B_2}, s_2s_1)$}
\end{picture}
\begin{picture}(140,80)(-40,-5)
\thicklines
\drawline[90](0,10)(80,10)
\drawline[90](0,30)(80,30)
\drawline[90](0,50)(80,50)
\drawline[90](0,70)(80,70)
  \put(0,10){\circle*{2.5}}
  \put(0,30){\circle*{2.5}}
  \put(0,50){\circle*{2.5}}
  \put(0,70){\circle*{2.5}}
  \put(80,10){\circle*{2.5}}
  \put(80,30){\circle*{2.5}}
  \put(80,50){\circle*{2.5}}
  \put(80,70){\circle*{2.5}}
\put(20,30){\circle*{2}}
\put(20,70){\circle*{2}}
\put(40,10){\circle*{2}}
\put(40,50){\circle*{2}}
\put(60,30){\circle*{2}}
\put(60,70){\circle*{2}}
  \put(-9,8){\scriptsize${1}$}
  \put(-9,28){\scriptsize${2}$}
  \put(-9,48){\scriptsize$\barred{2}$}
  \put(-9,68){\scriptsize${\barred{1}}$}
  \put(84,8){\scriptsize${1}$}
  \put(84,28){\scriptsize${2}$}
  \put(84,48){\scriptsize${\barred{2}}$}
  \put(84,68){\scriptsize${\barred{1}}$}
  \put(87,60){\oval(7,32)}
  \put(-7,60){\oval(7,32)}
\thicklines
\light{
\drawline[90](80,50)(60,30)
\drawline[90](60,70)(40,50)
\drawline[90](60,30)(40,10)
  \put(50,55){\scriptsize$t_1$}
  \put(50,15){\scriptsize$t_1$}
  \put(77,38){\scriptsize$t_2^2$}
}
\dark{
\drawline[90](40,50)(20,70)
\drawline[90](40,10)(20,30)
\drawline[90](20,30)(0,50)
}
\put(16,-8){$\Gamma_{\mathrm{S}} (\mathsf{B_2}, s_2s_1)$}
\end{picture}
\end{center}

\begin{center}
\begin{picture}(70,60)(5,-5)
\thinlines
\dottedline{2}(0,0)(40,0)
\dottedline{2}(0,10)(40,10)
\dottedline{2}(0,20)(40,20)
\dottedline{2}(0,30)(40,30)
\dottedline{2}(0,40)(40,40)
  \put(0,30){\circle*{2}}
  \put(40,30){\circle*{2}}
  \put(-8,28){\scriptsize$\barred{2}$}
  \put(42,28){\scriptsize$\barred{2}$}
\dottedline{2}(40,30)(30,20)
\dottedline{2}(40,30)(30,10)
\dottedline{2}(40,20)(30,10)
\dottedline{2}(30,40)(20,30)
\dottedline{2}(30,10)(20,0)
\dottedline{2}(20,30)(10,40)
\dottedline{2}(20,0)(10,10)
\dottedline{2}(10,10)(0,20)
\dottedline{2}(10,10)(0,30)
\dottedline{2}(10,20)(0,30)
\thicklines
\drawline[90](0,30)(40,30)
\put(5,-12){$\mathrm{wt}: 1$}
\end{picture}
\begin{picture}(70,60)(5,-5)
\thinlines
\dottedline{2}(0,0)(40,0)
\dottedline{2}(0,10)(40,10)
\dottedline{2}(0,20)(40,20)
\dottedline{2}(0,30)(40,30)
\dottedline{2}(0,40)(40,40)
  \put(0,30){\circle*{2}}
  \put(40,30){\circle*{2}}
      \put(-8,28){\scriptsize$\barred{2}$}
  \put(42,28){\scriptsize$\barred{2}$}
\dottedline{2}(40,30)(30,20)
\dottedline{2}(40,30)(30,10)
\dottedline{2}(40,20)(30,10)
\dottedline{2}(30,40)(20,30)
\dottedline{2}(30,10)(20,0)
\dottedline{2}(20,30)(10,40)
\dottedline{2}(20,0)(10,10)
\dottedline{2}(10,10)(0,20)
\dottedline{2}(10,10)(0,30)
\dottedline{2}(10,20)(0,30)
\thicklines
\drawline[90](40,30)(30,20)(10,20)(0,30)
\put(5,-12){$\mathrm{wt}: 2t_2$}
\end{picture}
\begin{picture}(60,30)(-20,-5)
\thinlines
\dottedline{2}(0,5)(40,5)
\dottedline{2}(0,15)(40,15)
\dottedline{2}(0,25)(40,25)
\dottedline{2}(0,35)(40,35)
  \put(0,25){\circle*{2}}
    \put(0,35){\circle*{2}}
  \put(40,25){\circle*{2}}
    \put(40,35){\circle*{2}}
  \put(-8,23){\scriptsize$\barred{2}$}
  \put(42,23){\scriptsize$\barred{2}$}
    \put(-8,33){\scriptsize$\barred{1}$}
  \put(42,33){\scriptsize$\barred{1}$}
\dottedline{2}(40,25)(20,5)(0,25)
\dottedline{2}(30,35)(20,25)(10,35)
\thicklines
\drawline[90](40,35)(0,35)
\drawline[90](0,25)(40,25)
\put(5,-12){$\mathrm{wt}: 1$}
\end{picture}
\begin{picture}(60,30)(-20,-5)
\thinlines
\dottedline{2}(0,5)(40,5)
\dottedline{2}(0,15)(40,15)
\dottedline{2}(0,25)(40,25)
\dottedline{2}(0,35)(40,35)
  \put(0,25){\circle*{2}}
    \put(0,35){\circle*{2}}
  \put(40,25){\circle*{2}}
    \put(40,35){\circle*{2}}
  \put(-8,23){\scriptsize$\barred{2}$}
  \put(42,23){\scriptsize$\barred{2}$}
    \put(-8,33){\scriptsize$\barred{1}$}
  \put(42,33){\scriptsize$\barred{1}$}
\dottedline{2}(40,25)(20,5)(0,25)
\dottedline{2}(30,35)(20,25)(10,35)
\thicklines
\drawline[90](40,35)(0,35)
\drawline[90](40,25)(30,15)(10,15)(0,25)
\put(5,-12){$\mathrm{wt}: t_2^2$}
\end{picture}
\end{center}

\begin{picture}(70,60)(35,-5)
\thinlines
\dottedline{2}(0,0)(40,0)
\dottedline{2}(0,10)(40,10)
\dottedline{2}(0,20)(40,20)
\dottedline{2}(0,30)(40,30)
\dottedline{2}(0,40)(40,40)
  \put(0,30){\circle*{2}}
  \put(40,30){\circle*{2}}
    \put(-8,28){\scriptsize$\barred{2}$}
  \put(42,28){\scriptsize$\barred{2}$}
\dottedline{2}(40,30)(30,20)
\dottedline{2}(40,30)(30,10)
\dottedline{2}(40,20)(30,10)
\dottedline{2}(30,40)(20,30)
\dottedline{2}(30,10)(20,0)
\dottedline{2}(20,30)(10,40)
\dottedline{2}(20,0)(10,10)
\dottedline{2}(10,10)(0,20)
\dottedline{2}(10,10)(0,30)
\dottedline{2}(10,20)(0,30)
\thicklines
\drawline[90](40,30)(30,10)(10,10)(0,30)
\put(5,-12){$\mathrm{wt}: t_2^2$}
\end{picture}
\begin{picture}(70,60)(35,-5)
\thinlines
\dottedline{2}(0,0)(40,0)
\dottedline{2}(0,10)(40,10)
\dottedline{2}(0,20)(40,20)
\dottedline{2}(0,30)(40,30)
\dottedline{2}(0,40)(40,40)
  \put(0,30){\circle*{2}}
  \put(40,30){\circle*{2}}
    \put(-8,28){\scriptsize$\barred{2}$}
  \put(42,28){\scriptsize$\barred{2}$}
\dottedline{2}(40,30)(30,20)
\dottedline{2}(40,30)(30,10)
\dottedline{2}(40,20)(30,10)
\dottedline{2}(30,40)(20,30)
\dottedline{2}(30,10)(20,0)
\dottedline{2}(20,30)(10,40)
\dottedline{2}(20,0)(10,10)
\dottedline{2}(10,10)(0,20)
\dottedline{2}(10,10)(0,30)
\dottedline{2}(10,20)(0,30)
\thicklines
\drawline[90](40,30)(30,10)(20,0)(10,10)(0,30)
\put(5,-12){$\mathrm{wt}: t_1t_2^2$}
\end{picture}
\begin{picture}(60,30)(-20,-5)
\thinlines
\dottedline{2}(0,5)(40,5)
\dottedline{2}(0,15)(40,15)
\dottedline{2}(0,25)(40,25)
\dottedline{2}(0,35)(40,35)
  \put(0,25){\circle*{2}}
    \put(0,35){\circle*{2}}
  \put(40,25){\circle*{2}}
    \put(40,35){\circle*{2}}
  \put(-8,23){\scriptsize$\barred{2}$}
  \put(42,23){\scriptsize$\barred{2}$}
    \put(-8,33){\scriptsize$\barred{1}$}
  \put(42,33){\scriptsize$\barred{1}$}
\dottedline{2}(40,25)(20,5)(0,25)
\dottedline{2}(30,35)(20,25)(10,35)
\thicklines
\drawline[90](40,35)(30,35)(20,25)(10,35)(0,35)
\drawline[90](40,25)(20,5)(0,25)
\put(5,-12){$\mathrm{wt}: (t_1t_2)^2$}
\end{picture}
\begin{center}
\begin{picture}(10,10)(0,0)
\dottedline(5,0)(5,230)
\end{picture}
\end{center}
\caption{$F_{c \omega_1} (t_1,t_2)$ and $F_{c^2 \omega_2}(t_1,t_2)$ in type $\mathsf{B_2}$ with $c=s_2s_1$.} \label{fig:bnexample}
\end{figure}
\end{example}

\begin{figure}[h]
\setlength{\unitlength}{1.2pt}
\begin{center}
\begin{picture}(63,120)(0,-15)
\thicklines
\put(18,0){\line(1,0){20}} \put(18,20){\line(1,0){20}} \put(18,40){\line(1,0){20}} \put(18,60){\line(1,0){20}}
\put(18,80){\line(1,0){20}} \put(18,100){\line(1,0){20}}
  \put(18,0){\circle*{2.5}}
  \put(18,20){\circle*{2.5}}
  \put(18,40){\circle*{2.5}}
  \put(18,60){\circle*{2.5}}
  \put(18,80){\circle*{2.5}}
  \put(18,100){\circle*{2.5}}
  \put(38,0){\circle*{2.5}}
  \put(38,20){\circle*{2.5}}
  \put(38,40){\circle*{2.5}}
  \put(38,60){\circle*{2.5}}
  \put(38,80){\circle*{2.5}}
  \put(38,100){\circle*{2.5}}
  \put(4,-2){\scriptsize${1}$}
  \put(4,18){\scriptsize${i}$}
  \put(0,38){\scriptsize${i+1}$}
  \put(0,58){\scriptsize${\overline{i+1}}$}
  \put(4,78){\scriptsize${\overline{i}}$}
  \put(4,98){\scriptsize${\overline{1}}$}
  \dottedline{3}(6,6)(6,13)
  \dottedline{3}(6,46)(6,53)
  \dottedline{3}(6,86)(6,93)
  \dottedline{3}(28,5)(28,15)
  \dottedline{3}(28,45)(28,55)
  \dottedline{3}(28,85)(28,95)
\dark{ \put(18,40){\line(1,-1){20}} \put(18,80){\line(1,-1){20}} }
  \put(22,26){\footnotesize$t$}
  \put(22,66){\footnotesize$t$}
  \put(17,-17){$x_{\barred{i}} (t)$}
\end{picture}
\begin{picture}(50,120)(0,-15)
\thicklines
\put(5,0){\line(1,0){20}} \put(5,20){\line(1,0){20}} \put(5,40){\line(1,0){20}} \put(5,60){\line(1,0){20}}
\put(5,80){\line(1,0){20}} \put(5,100){\line(1,0){20}}
  \put(5,0){\circle*{2.5}}
  \put(5,20){\circle*{2.5}}
  \put(5,40){\circle*{2.5}}
  \put(5,60){\circle*{2.5}}
  \put(5,80){\circle*{2.5}}
  \put(5,100){\circle*{2.5}}
  \put(25,0){\circle*{2.5}}
  \put(25,20){\circle*{2.5}}
  \put(25,40){\circle*{2.5}}
  \put(25,60){\circle*{2.5}}
  \put(25,80){\circle*{2.5}}
  \put(25,100){\circle*{2.5}}
  \dottedline{3}(15,5)(15,15)
  \dottedline{3}(15,45)(15,55)
  \dottedline{3}(15,85)(15,95)
\light{ \put(5,20){\line(1,1){20}} \put(5,60){\line(1,1){20}} }
  \put(18,26){\footnotesize$t$}
  \put(18,66){\footnotesize$t$}
  \put(4,-17){$x_i (t)$}
\end{picture}
\begin{picture}(50,120)(0,-15)
\thicklines
\put(5,0){\line(1,0){20}}  \put(5,40){\line(1,0){20}} \put(5,60){\line(1,0){20}}
 \put(5,100){\line(1,0){20}}
  \put(5,0){\circle*{2.5}}
  \put(5,40){\circle*{2.5}}
  \put(5,60){\circle*{2.5}}
  \put(5,100){\circle*{2.5}}
  \put(25,0){\circle*{2.5}}
  \put(25,40){\circle*{2.5}}
  \put(25,60){\circle*{2.5}}
  \put(25,100){\circle*{2.5}}
  \dottedline{3}(15,5)(15,35)
  \dottedline{3}(15,65)(15,95)
\dark{ \put(5,60){\line(1,-1){20}}}
  \put(11,44){\footnotesize$t$}
  \put(4,-17){$x_{\barred{n}} (t)$}
\end{picture}
\begin{picture}(50,120)(0,-15)
\thicklines
\put(5,0){\line(1,0){20}}  \put(5,40){\line(1,0){20}} \put(5,60){\line(1,0){20}}
 \put(5,100){\line(1,0){20}}
  \put(5,0){\circle*{2.5}}
  \put(5,40){\circle*{2.5}}
  \put(5,60){\circle*{2.5}}
  \put(5,100){\circle*{2.5}}
  \put(25,0){\circle*{2.5}}
  \put(25,40){\circle*{2.5}}
  \put(25,60){\circle*{2.5}}
  \put(25,100){\circle*{2.5}}
  \put(35,-2){\scriptsize${1}$}
  \put(35,38){\scriptsize${n}$}
  \put(35,58){\scriptsize${\overline{n}}$}
  \put(35,98){\scriptsize${\overline{1}}$}
  \dottedline{3}(15,5)(15,35)
  \dottedline{3}(15,65)(15,95)
  \dottedline{3}(36,8)(36,35)
  \dottedline{3}(36,68)(36,95)
\light{ \put(5,40){\line(1,1){20}}  }
  \put(18,45){\footnotesize$t$}
  \put(4,-17){$x_n (t)$}
\end{picture}
\end{center}
\caption{Elementary chips in $\Gamma (\mathsf{C_n}, c)$ ($i = 1, \ldots, n-1$).} \label{fig:cnchips}
\end{figure}
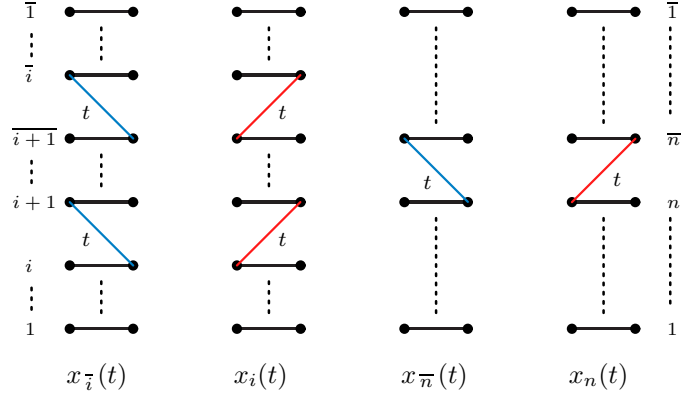

\textbf{Type $\mathsf{C_n}$} ($n \geq 2$) :
For each $i \in [1, n] \cup [\barred{1}, \barred{n}]$, the elementary chip corresponding to $x_{i}(t)$ is shown in Figure~\ref{fig:cnchips}. The vertices for each elementary chip consist of all the endpoints of the $2n$ horizontal edges. We number the horizontal levels from bottom to top by $1, \ldots, n, \barred{n}, \ldots, \barred{1}$. All the edges are oriented from right to left with weights shown in the figure, all unlabeled edges have weight 1.

The directed graph $\Gamma (\mathsf{C_n}, c)$ associated with $c = s_{i_1} \cdots s_{i_n}$ is constructed as a concatenation of elementary chips $x_{\barred{i_i}}(1), \ldots, x_{\barred{i_n}}(1), x_{i_n}(t_{i_n}), \ldots, x_{i_1}(t_{i_1})$ (in this order). We number the $2n$ sources and the $2n$ sinks of the graph $\Gamma (\mathsf{C_n}, c)$ bottom-to-top in the order $1, \ldots, n, \barred{n}, \ldots, \barred{1}$.
The definition of the weight of paths is the same as before.

The Weyl group of type $\mathsf{C_n}$ acts on the index set $[1, n] \cup [\barred{1}, \barred{n}]$ in the same way as the Weyl group of type $\mathsf{B_n}$. We then have the following proposition for computing the $F$-polynomials of type $\mathsf{C_n}$.

\begin{prop}
\label{pr:genminorsC} For $k \in [1,n]$ ,the $F$-polynomials $F_{c^m \omega_k}(t_1, \dots, t_n)$ equals the sum of weights of all collections of vertex-disjoint paths in $\Gamma (\mathsf{C_n}, c)$ with the sources and sinks labeled by~$c^m \cdot [1, k]$.
\end{prop}

\section{Proofs of the main results}
\label{sec:proofs}

We start by briefly recalling the definition of generalized minors; more details can be found in \cite{fz-double}.

Let $\lg$ be a complex semisimple Lie algebra of rank~$n$ with the \emph{Cartan decomposition} $\lg = \nn_- \oplus \hh \oplus \nn$.
Let $e_i, h_i, f_i\,$, for $i \in [1, n]$, be the standard generators of $\lg$.
The \emph{simple roots} $\alpha_i$ ($i \in [1, n]$) form a basis in the dual space $\hh^*$ such that  $[h, e_i] = \alpha_i (h) e_i$, and $[h,f_i] = - \alpha_i (h) f_i$
for any $h \in \hh$ and $i \in [1, n]$. The structure of $\lg$ is uniquely determined by the
\emph{Cartan matrix} $A = (a_{i,j})$ given by $a_{i,j} = \alpha_j (h_i)$.

Let $G$ be a simply connected complex Lie group with the Lie algebra $\lg$. For $i \in [1, n]$, let $\varphi_i : \SL_2 \rightarrow G$ denote the canonical embedding corresponding to the simple root $\alpha_i$. For $i \in [1, n]$ and $t \in \CC$, we write
\beal\label{eq:x,y}
x_i (t) = \varphi_i \mat{1}{t}{0}{1} = \exp(t e_i)\,, \quad x_{\barred{i}} (t) = \varphi_i \mat{1}{0}{t}{1} = \exp(t f_i)\,.
\eea
We also set
$$t^{h_i} = \varphi_i \mat{t}{0}{0}{t^{-1}}$$
for any $i \in [1, n]$ and any $t \neq 0$.
Let $N$ (resp., $N_-$) be the maximal unipotent subgroup of $G$ generated by all $x_i (t)$ (resp. $x_{\barred{i}} (t)$) with Lie algebra $\nn$ (resp., $\nn_-$), and $H$ be the maximal torus in $G$ with the Lie algebra $\hh$.

The \emph{Weyl group} $W$ of $G$ is defined to be the group of linear transformations of the root space $\hh^*$ generated by
the \emph{simple reflections} $s_1, \ldots, s_n$, whose action on $\hh^*$ is given by $s_i (\gamma) = \gamma - \gamma (h_i) \alpha_i$ for $\gamma \in \hh^*$.

A \emph{reduced word} for $w \in W$ is a sequence of indices
$(i_1, \ldots, i_m)$ of shortest possible length~$m$
such that $w = s_{i_1} \cdots s_{i_m}$.
The number~$m$ is denoted by $\l(w)$ and
is called the \emph{length} of~$w$. The group $W$ possesses a unique element $\wnot$ of maximal length.

The Weyl group $W$ is naturally identified with ${\rm Norm}_G (H)/H$ by sending each simple reflection
$s_i$ to the coset $\barred{s_i} H$, where the representative $\barred{s_i} \in {\rm Norm}_G (H)$ is defined by
\beal
\label{eq:s_i-via-x_i}
\barred{s_i} = \varphi_i \mat{0}{-1}{1}{0} =x_i(-1)x_{\barred{i}}(1) x_i(-1) \,.
\eea
The elements $\barred{s_i}$ satisfy the braid relations in~$W$; thus the representative $\barred{w}$ can be unambiguously defined for any
$w \in W$ by requiring that $\barred{uv} = \barred{u} \cdot \barred{v}$ whenever $\l (uv) = \l (u) + \l (v)$.

The \emph{weight lattice} $P$ is the set of all weights $\gamma \in \hh^*$
such that $\gamma (h_i) \in \ZZ$ for all $i$.
The group $P$ has a $\ZZ$-basis formed by the \emph{fundamental weights}
$\omega_1, \ldots, \omega_n$ defined by $\omega_i (h_j) = \delta_{ij}$.
With some abuse of notation, we identify the weight lattice $P$ in $\hh^*$
with the group of rational multiplicative characters of~$H$,
here written in the exponential notation:
a weight $\gamma \in P$ acts by $a \mapsto a^\gamma$.
Under this identification, the fundamental weights $\omega_1, \ldots, \omega_n$ act in $H$ by
$(t^{h_j})^{\omega_i} = t^{\delta_{ij}}$.
Recall that the set $G_0 = N_-HN$ of elements $x\in G$ that
have Gaussian decomposition is open and dense in~$G$. This (unique) decomposition of
$x \in N_-HN$ will be written as $x = [x]_- [x]_0 [x]_+$\,.

We now define the \emph{generalized minors} introduced in \cite{fz-double}.
For $u,v \in W$ and $k \in [1, n]$, the generalized minor
$\Delta_{u \omega_k, v \omega_k}$
is the regular function on $G$ whose restriction to the open set
${\barred{u}} G_0 {\barred{v}}^{-1}$ is given by
\begin{equation}
\label{eq:Delta-general}
\Delta_{u \omega_k, v \omega_k} (x) =
(\left[{\barred{u}}^{-1}
   x {\barred{v}}\right]_0)^{\omega_k}\,.
\end{equation}
As shown in \cite{fz-double}, $\Delta_{u \omega_k, v \omega_k}$ depends on
the weights $u \omega_k$ and $v \omega_k$ alone, not on the particular
choice of $u$ and~$v$.
Let $V_{\omega_k}$ be the fundamental representation of $G$ and $\vv$ be a highest weight vector with highest weight $\omega_k$, then $\barred{u} \mathbf{v}$ and $\barred{v} \mathbf{v}$ are two vectors with weights $u \omega_k$ and $v \omega_k$, respectively.
From the definition of generalized minors, it is not hard to see that $\Delta_{u \omega_k, v \omega_k} (x)$ is the coefficient of $\barred{u} \mathbf{v}$ in the expression of $x \cdot \barred{v} \mathbf{v}$ (the action of the group element $x$ on $V_{\omega_k}$) in terms of a weight basis containing both $\barred{u} \mathbf{v}$ and $\barred{v} \mathbf{v}$.

Let $c = s_{i_1} \cdots s_{i_n}$ be a Coxeter element and $x_c = x_{\barred{i_i}}(1) \cdots x_{\barred{i_n}}(1) x_{i_n}(t_{i_n}) \cdots  x_{i_1}(t_{i_1})$. We compute the generalized minors of the form $\Delta_{c^m \omega_k, c^m \omega_k} (x_c)$ (hence the $F$-polynomials) by explicitly computing the action by $x_{c}$ on each fundamental representation $V_{\omega_k}$; recall that $x_i(t)$ and $x_{\barred{i}}(t)$ act in every finite-dimensional  representation of $\lg$ by
\begin{equation}
\label{eq:xiinU} x_i (t) = \sum_{n \geq 0} \frac{t^n}{n !} \, e_i^{n} \quad \mbox{and} \quad x_{\barred{i}} (t) = \sum_{n \geq 0}
\frac{t^n}{n !} \, f_i^{n}\,.
\end{equation}

For the type $\mathsf{A_n}$ case (i.e., when $G = \SL_{n+1}$), the generalized minors specialize to the ordinary minors as follows.
The Weyl group $W$ is identified with the symmetric group~$S_{n+1}$, and
$V_{\omega_k} = \bigwedge^k \CC^{n+1}$, the $k$-th exterior power of the standard representation.
All the weights of $V_{\omega_k}$ are extremal, and are in bijection with the $k$-subsets of $[1,n+1]$, so that~$W$
acts on them in a natural way, and $\omega_k$ corresponds to~$[1,k]$.
If~$\gamma$ and~$\delta$ correspond to $k$-subsets $I$ and $J$, respectively,
then $\Delta_{\gamma,\delta} = \Delta_{I,J}$ is the minor with the row set~$I$ and the column set~$J$.

Note that the directed graph $\Gamma (\mathsf{A_n}, c)$ provides a combinatorial model for the action of $x_c$ in each $\bigwedge^k \CC^{n+1}$,
in the sense that each of the elementary chips corresponding to $x_{i} (t)$ captures the action of $x_{i}(t)$ on the fundamental representations. For example, let $G = \mathrm{SL}_3$ and $V$ be its standard representation with basis $\vv_1, \vv_2, \vv_3$. We have $x_1(t) \cdot (\vv_2 \wedge \vv_3) = \vv_2 \wedge \vv_3 + t \, \vv_1 \wedge \vv_3$, where the coefficient 1 (resp., $t$) of $\vv_2 \wedge \vv_3$ (resp., $\vv_1 \wedge \vv_3$) is the product of the weights of the edges connecting the sources labeled $\{2, 3\}$ and the sinks labeled by $\{2, 3\}$ (resp., $\{1, 3\}$).
The directed graphs $\Gamma (\mathsf{D_n}, c)$, $\Gamma (\mathsf{B_n}, c)$, $\Gamma_{\mathrm{S}} (\mathsf{B_n}, c)$ and $\Gamma (\mathsf{C_n}, c)$ are designed and constructed to serve the same purpose, that is to capture the action of $x_c$ on the fundamental representations. This will become clear after we recall the Lie algebra action on the corresponding fundamental representations in each of the classical types (c.f. \cite{hk}).


\textbf{Proof of Proposition~\ref{pr:genminorsD}:} Let $\lg$ be the simple Lie algebra of type $\mathsf{D_n}$ for $n \geq 4$, that is, the even special orthogonal Lie algebra $\so_{2n}$.  Then the action of generators in the standard $2n$-dimensional representation $V$ with respect to the standard basis $\vv_1, \ldots, \vv_n, \vv_{\barred{n}}, \ldots, \vv_{\barred{1}}$ can be written as:
\beal \label{soeaction}
e_i \cdot \vv_j =
\bmcase
    \vv_i\,, & \mbox{ if } i\neq n \mbox{ and } j = i+1\,;\\
    \vv_{\barred{i+1}}\,, & \mbox{ if } i \neq n \mbox{ and } j = \barred{i}\,;\\
    \vv_n\,, & \mbox{ if } i = n \mbox{ and } j = \barred{n-1}\,;\\
    \vv_{n-1}\,, & \mbox{ if } i = n \mbox{ and } j = \barred{n}\,;\\
    0\,, & \mbox{ otherwise}, \\
\emcase \\[.5in]
f_i \cdot \vv_j =
\bmcase
    \vv_{i+1}\,, & \mbox{ if } i \neq n \mbox{ and } j = i\,; \\
    \vv_{\barred{i}}\,, & \mbox{ if } i \neq n \mbox{ and } j = \barred{i+1}\,;\\
    \vv_{\barred{n}}\,, & \mbox{ if } i = n \mbox{ and } j = n-1\,;\\
    \vv_{\barred{n-1}}\,, & \mbox{ if } i = n \mbox{ and } j = n\,;\\
    0\,, & \mbox{ otherwise}. \\
\emcase
\eea
for $i \in [1, n]$ and $j \in [1, n] \cup [\barred{1}, \barred{n}]$.

The group elements $x_i(t)$ and $x_\barred{i}(t)$ act as $I + t e_i$ and $I + t f_i$ respectively on $V$.
%
%
We associate each vertex of $\Gamma (\mathsf{D_n}, c)$ on the horizontal level $j$ with the basis vector $\vv_j \in V$ for $j \in [1,n] \cup [\barred{1}, \barred{n}]$, then the action of $x_i (t)$ and $x_{\barred{i}}(t)$ on $V$ can be read from the corresponding elementary chips.
For instance, the fragment shown in Figure~\ref{fig:partofxi} expresses the action $x_{i}(t) \cdot \vv_{\light{\barred{i}}} = 1 \vv_{\dark{\barred{i}}} + t \vv_{\dgreen{\barred{i+1}}}$ for $i \neq n$. Note that this fragment is part of the elementary chip corresponding to $x_{i}(t)$ for $i \neq n$.
\begin{figure}[h]
\setlength{\unitlength}{1.2pt}
\begin{center}
\begin{picture}(20,30)(0,0)
\drawline[90](0,0)(20,0)
\drawline[120](0,20)(20,20)
\drawline[120](20,20)(0,0)
\thicklines
  \put(0,0){\circle*{2}}
  \put(0,20){\circle*{2}}
  \put(20,0){\circle*{1.5}}
  \put(20,20){\circle*{2.5}}
  \put(10,23){\footnotesize$1$}
  \put(11,5){\footnotesize$t$}
  \dark{\put(-10,18){\scriptsize$\barred{i}$}}
  \dgreen{\put(-20,-2){\scriptsize$\barred{i+1}$}}
  \light{\put(24,18){\scriptsize$\barred{i}$}}
\end{picture}
\end{center}
\caption{} \label{fig:partofxi}
\end{figure}
Therefore the graph $\Gamma (\mathsf{D_n}, c)$ (constructed by concatenation of the elementary chips) provides a combinatorial model for the action of $x_c$ on $V$, that is, the coefficient of $\vv_j$ in the expression of $x_c \cdot \vv_i$ is equal to the sum of the weights of all paths in $\Gamma (\mathsf{D_n}, c)$ with source labeled by $i$ and sink labeled by $j$.

This observation can be generalized to the exterior powers of $V$ and used to compute the generalized minors. Recall that in the type $\mathsf{D_n}$ case, the fundamental representation $V_{\omega_k}$ for $k = 1, \ldots, n-2$ is realized as $\bigwedge^k V$ with the highest weight vector $\vv_1 \wedge \cdots \wedge \vv_{k}$. Each extremal weight $u\omega_k$ of $V_{\omega_k}$ corresponds to a $k$-subset $u\cdot [1, k]$ in $[1, n] \cup [\barred{1}, \barred{n}]$. Note that $i$ and $\barred{i}$ do not appear simultaneously in $u \cdot [1, k]$ for any $i \in [1, n]$ and $u \in W$. We define a linear ordering on the index set $[1, n] \cup [\barred{1}, \barred{n}]$ by $1 < \cdots < n < \barred{n} < \cdots < 1$. Let $I = \{i_1 < \cdots < i_k\}$ be a $k$-subset in $[1, n] \cup [\barred{1}, \barred{n}]$ corresponding to an extremal weight $\gamma$, and define a basis vector $\vv_{I} = \vv_{i_1} \wedge \cdots \wedge \vv_{i_k}$ in $\bigwedge^k V$. Then the principal minor $\Delta_{\gamma, \gamma}(x_c)$ equals to the coefficient of $\vv_{I}$ in the expression of $x_c \cdot \vv_I$ (in terms of the standard basis in $\bigwedge^k V$). It can be computed as follows:

\begin{quote}
\label{eq:pminpaths}
$\Delta_{\gamma, \gamma}(x_c)$ equals the sum of \textbf{signed-weights} of all collections of vertex-disjoint paths in $\Gamma (\mathsf{D_n}, c)$ with sources and sinks labeled by $I$.
\end{quote}

The requirement of the paths to be vertex-disjoint is because $\vv \wedge \vv = 0$ for any $\vv \in V$.

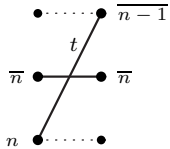
\begin{figure}[h]
\setlength{\unitlength}{1.2pt}
\begin{center}
\begin{picture}(20,50)(0,0)
  \dottedline{3}(0,0)(20,0)
  \dottedline{3}(0,40)(20,40)
\thicklines
\drawline[90](0,0)(20,40)
\drawline[90](0,20)(20,20)
  \put(0,0){\circle*{2.5}}
  \put(0,20){\circle*{2.5}}
  \put(20,20){\circle*{2.5}}
  \put(20,40){\circle*{2.5}}
  \put(0,40){\circle*{2}}
  \put(20,0){\circle*{2}}
  \put(10,28){\footnotesize$t$}
  \put(-10,-2){\scriptsize$n$}
  \put(-10,18){\scriptsize$\barred{n}$}
  \put(24,18){\scriptsize$\barred{n}$}
  \put(24,38){\scriptsize$\barred{n-1}$}
\end{picture}
\end{center}
\caption{crossing happen in the chip $x_{n}(t)$ on level $\barred{n}$.} \label{fig:crossing}
\end{figure}

To define the signed-weight of a family of paths, we first recall that in the chips corresponding to $x_{\barred{n}}(t)$ and $x_n(t)$, the intersections of the diagonal edges and the horizontal edges in the middle of each horizontal edge are \emph{not} vertices. Hence two vertex-disjoint paths in $\Gamma (\mathsf{D_n}, c)$ can \emph{cross} each other at such points (see Figure~\ref{fig:dnchipsnsp}).
One crossing of this kind is shown in Figure~\ref{fig:crossing} and the two paths crossing each other are depicted by thick lines  (Note that there are four kinds of crossings in $\Gamma (\mathsf{D_n}, c)$, see Figure~\ref{fig:p2}).
It represents that the expression of $x_{n}(t) \cdot (\vv_{\barred{n}} \wedge \vv_{\barred{n-1}})$ contains the term $t \vv_{\barred{n}} \wedge \vv_{n} = -t \vv_{n} \wedge \vv_{\barred{n}}$. This negative coefficient leads to the definition of signed-weight.
We define the \emph{signed-weight} of a family of paths to be the weight of the family of paths if there are an~even number of such crossings in the paths, and to be the \emph{negative} of the weight of the family of paths if there are an~odd number of such crossings.

The crossing can only happen in the elementary chips corresponding to $x_{\barred{n}}(t)$ and $x_{n}(t)$ and this two chips appear in $\Gamma (\mathsf{D_n}, c)$ exactly once, therefore at most two crossings can appear in a family of vertex-disjoint paths in $\Gamma (\mathsf{D_n}, c)$. Also when the sources and the sinks are labeled by the same index set, we always have a family of vertex-disjoint paths consisting of the horizontal levels connecting the sources and the sinks. This family has (signed-)weight $1$.
Hence,
to prove part $(1)$ of Propositions~\ref{pr:genminorsD}, it is enough to show that there does not exist a family of vertex-disjoint paths in $\Gamma (\mathsf{D_n}, c)$ with the sources and sinks labeled by~$c^m \cdot [1, k]$ for $k \in [2, n-2]$ such that there is exactly one crossing among its paths.

Suppose for the sake of contradiction that such a family of paths exists.
We use $x_i$ to represent the corresponding element chip if there is no danger of confusion. All families of paths are assumed to be vertex-disjoint. We first consider the case that the (only one) crossing appears in $x_n$ on the level $\barred{n}$.

\begin{figure}[ht]
\setlength{\unitlength}{1.2pt}
\begin{center}
\begin{picture}(140,70)(0,-15)
\thinlines
\dottedline{2}(0,0)(140,0)
\dottedline{2}(0,10)(140,10)
\dottedline{2}(0,20)(140,20)
\dottedline{2}(0,30)(140,30)
  \put(0,20){\circle*{2.5}}
  \put(0,50){\circle*{2.5}}
  \put(140,20){\circle*{2.5}}
  \put(140,50){\circle*{2.5}}
    \put(20,20){\circle*{2}}
    \put(20,30){\circle*{2}}
    \put(30,10){\circle*{2}}
    \put(30,0){\circle*{2}}
    \put(50,10){\circle*{2}}
    \put(50,30){\circle*{2}}
    \put(60,0){\circle*{2}}
    \put(60,20){\circle*{2}}
    \put(80,0){\circle*{2}}
    \put(80,20){\circle*{2}}
    \put(90,10){\circle*{2}}
    \put(90,30){\circle*{2}}
    \put(110,0){\circle*{2}}
    \put(110,10){\circle*{2}}
    \put(120,20){\circle*{2}}
    \put(120,30){\circle*{2}}
\put(-10,18){\scriptsize${\barred{n}}$}
\put(-20,28){\scriptsize${\barred{n-1}}$}
\put(-20,48){\scriptsize${\barred{n-i}}$}
\put(-10,8){\scriptsize${n}$}
\put(-20,-2){\scriptsize$n-1$}
\put(144,18){\scriptsize${\barred{n}}$}
\put(144,28){\scriptsize${\barred{n-1}}$}
\put(144,48){\scriptsize${\barred{n-i}}$}
\put(144,8){\scriptsize${n}$}
\put(144,-2){\scriptsize${n-1}$}
\put(48,-10){\scriptsize$x_{\barred{n-1}}$}
\put(20,-10){\scriptsize$x_{\barred{n}}$}
\put(110,-10){\scriptsize$x_{n}$}
\put(80,-10){\scriptsize$x_{n-1}$}
\dottedline{2}(120,30)(110,10)
\dottedline{2}(120,20)(110,0)
\dottedline{2}(90,30)(80,20)
\dottedline{2}(90,10)(80,0)
\dottedline{2}(60,20)(50,30)
\dottedline{2}(60,0)(50,10)
\dottedline{2}(30,10)(20,30)
\dottedline{2}(30,0)(20,20)
\thicklines
\light{
\drawline[90](140,20)(60,20)(50,30)(20,30)
\dottedline{3}(20,30)(0,50)
\put(129,15){\scriptsize$p_2$}
}
\Thicklines
\dark{
\dottedline{3}(140,50)(120,30)
\drawline[90](120,30)(110,10)(90,10)(80,0)
\dottedline{3}(80,0)(60,0)
\dottedline{3}(80,0)(70,-10)(60,0)
\drawline[90](60,0)(30,0)(20,20)(0,20)
\put(126,45){\scriptsize$p_1$}
}
\end{picture}
\end{center}
\caption{ }
\label{fig:p1}
\end{figure}

Let $\mathcal{P}$ be a family of paths and $p_1$, $p_2$ be two paths in $\mathcal{P}$ crossing each other in $x_n$ on the level $\barred{n}$. Let $p_1$ be the path that passes through the (upper) diagonal edge of $x_n$, and let $p_2$ be the path that passes through the horizontal edge of $x_n$ on level $\barred{n}$. We claim that the paths $p_1$ and $p_2$ appear partially as illustrated in Figure~\ref{fig:p1}. Recall that all paths travel from right to left. In Figure~\ref{fig:p1}, if $i \geq 2$ then $p_1$ must stay on level $n-1$ between $x_{\barred{n-1}}$ and $x_{n-1}$. When $i = 1$, it is possible for $p_1$ to go down after $x_{n-1}$, then go up to the level $n-1$ before $x_{\barred{n-1}}$.

Since $p_1$ and $p_2$ are the only paths crossing each other, it is easy to see that the label of the source of $p_1$ (resp., $p_2$) will be the label of the sink of $p_2$ (resp., $p_1$). Also the source of $p_1$ is ``higher'' than the source of $p_2$ (i.e, the label of the source of $p_1$ is bigger than the label of the source of $p_2$ in the linear order on $[1, n] \cup [\barred{1}, \barred{n}]$ defined before), since for $i \in [1, n-1]$, all the edges of the chip $x_i$ either keep the same horizontal level or bring the level down by 1. Each $x_i$ appears exactly once in $\Gamma (\mathsf{D_n}, c)$, hence the labels of the sources of $\{p_1, p_2\}$ must be $\{\barred{n-i}, \barred{n}\}$ or $\{\barred{n-1-i}, \barred{n-1}\}$ for some $i \geq 1$. To see the later case cannot happen, we assume that the labels of the sources of $\{p_1, p_2\}$ are $\{\barred{n-i}, \barred{n-1}\}$ for some $i \geq 2$: In this case, the chips $x_{\barred{n-1}}, x_{\barred{n-2}}, x_{\barred{n}}, x_{n}, x_{n-2}, x_{n-1}$ must appear in $\Gamma (\mathsf{D_n}, c)$ in this order. See Figure~\ref{fig:pf}.

\begin{figure}[ht]
\setlength{\unitlength}{1.2pt}
\begin{center}
\begin{picture}(200,90)(0,-15)
\thinlines
\dottedline{2}(0,0)(200,0)
\dottedline{2}(0,10)(200,10)
\dottedline{2}(0,20)(200,20)
\dottedline{2}(0,30)(200,30)
\dottedline{2}(0,40)(200,40)
\dottedline{2}(0,50)(200,50)
  \put(0,40){\circle*{2.5}}
  \put(0,70){\circle*{2.5}}
  \put(200,40){\circle*{2.5}}
  \put(200,70){\circle*{2.5}}
    \put(20,20){\circle*{2}}
    \put(20,40){\circle*{2}}
    \put(30,10){\circle*{2}}
    \put(30,30){\circle*{2}}
    \put(50,10){\circle*{2}}
    \put(50,50){\circle*{2}}
    \put(60,0){\circle*{2}}
    \put(60,40){\circle*{2}}
    \put(80,30){\circle*{2}}
    \put(80,40){\circle*{2}}
    \put(90,10){\circle*{2}}
    \put(90,20){\circle*{2}}
    \put(110,10){\circle*{2}}
    \put(110,20){\circle*{2}}
    \put(120,30){\circle*{2}}
    \put(120,40){\circle*{2}}
    \put(140,0){\circle*{2}}
    \put(140,40){\circle*{2}}
    \put(150,10){\circle*{2}}
    \put(150,50){\circle*{2}}
    \put(170,10){\circle*{2}}
    \put(170,30){\circle*{2}}
    \put(180,20){\circle*{2}}
    \put(180,40){\circle*{2}}
\put(-10,28){\scriptsize${\barred{n}}$}
\put(-20,38){\scriptsize${\barred{n-1}}$}
\put(-20,48){\scriptsize${\barred{n-2}}$}
\put(-20,68){\scriptsize${\barred{n-i}}$}
\put(-10,18){\scriptsize${n}$}
\put(-20,8){\scriptsize$n-1$}
\put(-20,-2){\scriptsize$n-2$}
\put(204,28){\scriptsize${\barred{n}}$}
\put(204,38){\scriptsize${\barred{n-1}}$}
\put(204,48){\scriptsize${\barred{n-2}}$}
\put(204,68){\scriptsize${\barred{n-i}}$}
\put(204,18){\scriptsize${n}$}
\put(204,8){\scriptsize$n-1$}
\put(204,-2){\scriptsize${n-2}$}
\put(16,-10){\scriptsize$x_{\barred{n-1}}$}
\put(46,-10){\scriptsize$x_{\barred{n-2}}$}
\put(80,-10){\scriptsize$x_{\barred{n}}$}
\put(110,-10){\scriptsize$x_{n}$}
\put(136,-10){\scriptsize$x_{n-2}$}
\put(166,-10){\scriptsize$x_{n-1}$}
\dottedline{2}(180,40)(170,30)
\dottedline{2}(180,20)(170,10)
\dottedline{2}(150,50)(140,40)
\dottedline{2}(150,10)(140,0)
\dottedline{2}(120,40)(110,20)
\dottedline{2}(120,30)(110,10)
\dottedline{2}(90,20)(80,40)
\dottedline{2}(90,10)(80,30)
\dottedline{2}(60,0)(50,10)
\dottedline{2}(60,40)(50,50)
\dottedline{2}(30,10)(20,20)
\dottedline{2}(30,30)(20,40)
\thicklines
\light{
\drawline[90](200,40)(180,40)(170,30)(80,30)
\put(190,35){\scriptsize$p_2$}
}
\Thicklines
\dark{
\dottedline{3}(200,70)(150,50)
\drawline[90](150,50)(140,40)(120,40)(110,20)
\put(182,68){\scriptsize$p_1$}
}
\end{picture}
\end{center}
\caption{ }
\label{fig:pf}
\end{figure}

Now consider the path $p_2$. It stays on the horizontal level $\barred{n}$ in the chip $x_{n}$. There is no edge connecting the horizontal levels $\barred{n}$ and $n$, and the chip $x_{\barred{n-2}}$ appears on the left of $x_{\barred{n}}$: therefore, we conclude that $p_2$ must stay on the horizontal level $\barred{n}$ when it arrives in the chip $x_{\barred{n}}$. Then it is clearly impossible for the path $p_2$ to reach its sink $\barred{n-i}$.
Therefore the sources of $\{p_1, p_2\}$ must be $\{\barred{n-i}, \barred{n}\}$ for some $i \geq 1$ and it is easy to see that the paths $p_1$ and $p_2$ appear (partially) as shown in Figure~\ref{fig:p1}.

By a similar argument, we obtain in Figure~\ref{fig:p2} all cases of paths that have exactly one crossing in $\Gamma (\mathsf{D_n}, c)$. In all cases, we denote $p_1$ and $p_2$ to be the two paths crossing each other, where $p_1$ is the path having higher source.

\begin{figure}[ht]
\setlength{\unitlength}{1.2pt}
\begin{center}
\begin{picture}(140,70)(10,-15)
\thinlines
\dottedline{2}(0,0)(140,0)
\dottedline{2}(0,10)(140,10)
\dottedline{2}(0,20)(140,20)
\dottedline{2}(0,30)(140,30)
  \put(0,20){\circle*{2.5}}
  \put(0,50){\circle*{2.5}}
  \put(140,20){\circle*{2.5}}
  \put(140,50){\circle*{2.5}}
    \put(20,20){\circle*{2}}
    \put(20,30){\circle*{2}}
    \put(30,10){\circle*{2}}
    \put(30,0){\circle*{2}}
    \put(50,10){\circle*{2}}
    \put(50,30){\circle*{2}}
    \put(60,0){\circle*{2}}
    \put(60,20){\circle*{2}}
    \put(80,0){\circle*{2}}
    \put(80,20){\circle*{2}}
    \put(90,10){\circle*{2}}
    \put(90,30){\circle*{2}}
    \put(110,0){\circle*{2}}
    \put(110,10){\circle*{2}}
    \put(120,20){\circle*{2}}
    \put(120,30){\circle*{2}}
\put(-10,18){\scriptsize${\barred{n}}$}
\put(-20,28){\scriptsize${\barred{n-1}}$}
\put(-20,48){\scriptsize${\barred{n-i}}$}
\put(-10,8){\scriptsize${n}$}
\put(-20,-2){\scriptsize$n-1$}
\put(48,-10){\scriptsize$x_{\barred{n-1}}$}
\put(20,-10){\scriptsize$x_{\barred{n}}$}
\put(110,-10){\scriptsize$x_{n}$}
\put(80,-10){\scriptsize$x_{n-1}$}
\dottedline{2}(120,30)(110,10)
\dottedline{2}(120,20)(110,0)
\dottedline{2}(90,30)(80,20)
\dottedline{2}(90,10)(80,0)
\dottedline{2}(60,20)(50,30)
\dottedline{2}(60,0)(50,10)
\dottedline{2}(30,10)(20,30)
\dottedline{2}(30,0)(20,20)
\thicklines
\light{
\drawline[90](140,20)(60,20)(50,30)(20,30)
\dottedline{3}(20,30)(0,50)
\put(130,15){\scriptsize$p_2$}
}
\Thicklines
\dark{
\dottedline{3}(140,50)(120,30)
\drawline[90](120,30)(110,10)(90,10)(80,0)
\dottedline{3}(80,0)(60,0)
\dottedline{3}(80,0)(70,-10)(60,0)
\drawline[90](60,0)(30,0)(20,20)(0,20)
\put(126,45){\scriptsize$p_1$}
}
\end{picture}
\begin{picture}(140,70)(-10,-15)
\thinlines
\dottedline{2}(0,0)(140,0)
\dottedline{2}(0,10)(140,10)
\dottedline{2}(0,20)(140,20)
\dottedline{2}(0,30)(140,30)
  \put(0,20){\circle*{2.5}}
  \put(0,50){\circle*{2.5}}
  \put(140,20){\circle*{2.5}}
  \put(140,50){\circle*{2.5}}
    \put(20,20){\circle*{2}}
    \put(20,30){\circle*{2}}
    \put(30,10){\circle*{2}}
    \put(30,0){\circle*{2}}
    \put(50,10){\circle*{2}}
    \put(50,30){\circle*{2}}
    \put(60,0){\circle*{2}}
    \put(60,20){\circle*{2}}
    \put(80,0){\circle*{2}}
    \put(80,20){\circle*{2}}
    \put(90,10){\circle*{2}}
    \put(90,30){\circle*{2}}
    \put(110,0){\circle*{2}}
    \put(110,10){\circle*{2}}
    \put(120,20){\circle*{2}}
    \put(120,30){\circle*{2}}
\put(144,18){\scriptsize${\barred{n}}$}
\put(144,28){\scriptsize${\barred{n-1}}$}
\put(144,48){\scriptsize${\barred{n-i}}$}
\put(144,8){\scriptsize${n}$}
\put(144,-2){\scriptsize${n-1}$}
\put(48,-10){\scriptsize$x_{\barred{n-1}}$}
\put(20,-10){\scriptsize$x_{\barred{n}}$}
\put(110,-10){\scriptsize$x_{n}$}
\put(80,-10){\scriptsize$x_{n-1}$}
\dottedline{2}(120,30)(110,10)
\dottedline{2}(120,20)(110,0)
\dottedline{2}(90,30)(80,20)
\dottedline{2}(90,10)(80,0)
\dottedline{2}(60,20)(50,30)
\dottedline{2}(60,0)(50,10)
\dottedline{2}(30,10)(20,30)
\dottedline{2}(30,0)(20,20)
\thicklines
\light{
\drawline[90](140,20)(120,20)(110,0)(80,0)
\dottedline{3}(80,0)(60,0)
\dottedline{3}(80,0)(70,-10)(60,0)
\drawline[90](60,0)(50,10)(30,10)(20,30)
\dottedline{3}(20,30)(0,50)
\put(130,15){\scriptsize$p_2$}
}
\Thicklines
\dark{
\dottedline{3}(140,50)(120,30)
\drawline[90](120,30)(90,30)(80,20)(0,20)
\put(126,45){\scriptsize$p_1$}
}
\end{picture}
\end{center}

\begin{center}
\begin{picture}(140,70)(10,-15)
\thinlines
\dottedline{2}(0,0)(140,0)
\dottedline{2}(0,10)(140,10)
\dottedline{2}(0,20)(140,20)
\dottedline{2}(0,30)(140,30)
  \put(0,10){\circle*{2.5}}
  \put(0,50){\circle*{2.5}}
  \put(140,10){\circle*{2.5}}
  \put(140,50){\circle*{2.5}}
    \put(20,10){\circle*{2}}
    \put(20,30){\circle*{2}}
    \put(30,20){\circle*{2}}
    \put(30,0){\circle*{2}}
    \put(50,20){\circle*{2}}
    \put(50,30){\circle*{2}}
    \put(60,0){\circle*{2}}
    \put(60,10){\circle*{2}}
    \put(80,0){\circle*{2}}
    \put(80,10){\circle*{2}}
    \put(90,20){\circle*{2}}
    \put(90,30){\circle*{2}}
    \put(110,0){\circle*{2}}
    \put(110,20){\circle*{2}}
    \put(120,10){\circle*{2}}
    \put(120,30){\circle*{2}}
\put(-10,18){\scriptsize${\barred{n}}$}
\put(-20,28){\scriptsize${\barred{n-1}}$}
\put(-20,48){\scriptsize${\barred{n-i}}$}
\put(-10,8){\scriptsize${n}$}
\put(-20,-2){\scriptsize$n-1$}
\put(50,-10){\scriptsize$x_{\barred{n}}$}
\put(20,-10){\scriptsize$x_{\barred{n-1}}$}
\put(110,-10){\scriptsize$x_{n-1}$}
\put(80,-10){\scriptsize$x_{n}$}
\dottedline{2}(120,30)(110,20)
\dottedline{2}(120,10)(110,0)
\dottedline{2}(90,30)(80,10)
\dottedline{2}(90,20)(80,0)
\dottedline{2}(60,10)(50,30)
\dottedline{2}(60,0)(50,20)
\dottedline{2}(30,20)(20,30)
\dottedline{2}(30,0)(20,10)
\thicklines
\light{
\drawline[90](140,10)(60,10)(50,30)(20,30)
\dottedline{3}(20,30)(0,50)
\put(130,13){\scriptsize$p_2$}
}
\Thicklines
\dark{
\dottedline{3}(140,50)(120,30)
\drawline[90](120,30)(110,20)(90,20)(80,0)
\dottedline{3}(80,0)(60,0)
\dottedline{3}(80,0)(70,-10)(60,0)
\drawline[90](60,0)(30,0)(20,10)(0,10)
\put(126,45){\scriptsize$p_1$}
}
\end{picture}
\begin{picture}(140,70)(-10,-15)
\thinlines
\dottedline{2}(0,0)(140,0)
\dottedline{2}(0,10)(140,10)
\dottedline{2}(0,20)(140,20)
\dottedline{2}(0,30)(140,30)
  \put(0,10){\circle*{2.5}}
  \put(0,50){\circle*{2.5}}
  \put(140,10){\circle*{2.5}}
  \put(140,50){\circle*{2.5}}
    \put(20,10){\circle*{2}}
    \put(20,30){\circle*{2}}
    \put(30,20){\circle*{2}}
    \put(30,0){\circle*{2}}
    \put(50,20){\circle*{2}}
    \put(50,30){\circle*{2}}
    \put(60,0){\circle*{2}}
    \put(60,10){\circle*{2}}
    \put(80,0){\circle*{2}}
    \put(80,10){\circle*{2}}
    \put(90,20){\circle*{2}}
    \put(90,30){\circle*{2}}
    \put(110,0){\circle*{2}}
    \put(110,20){\circle*{2}}
    \put(120,10){\circle*{2}}
    \put(120,30){\circle*{2}}
\put(144,18){\scriptsize${\barred{n}}$}
\put(144,28){\scriptsize${\barred{n-1}}$}
\put(144,48){\scriptsize${\barred{n-i}}$}
\put(144,8){\scriptsize${n}$}
\put(144,-2){\scriptsize${n-1}$}
\put(50,-10){\scriptsize$x_{\barred{n}}$}
\put(20,-10){\scriptsize$x_{\barred{n-1}}$}
\put(110,-10){\scriptsize$x_{n-1}$}
\put(80,-10){\scriptsize$x_{n}$}
\dottedline{2}(120,30)(110,20)
\dottedline{2}(120,10)(110,0)
\dottedline{2}(90,30)(80,10)
\dottedline{2}(90,20)(80,0)
\dottedline{2}(60,10)(50,30)
\dottedline{2}(60,0)(50,20)
\dottedline{2}(30,20)(20,30)
\dottedline{2}(30,0)(20,10)
\thicklines
\light{
\drawline[90](140,10)(120,10)(110,0)(80,0)
\dottedline{3}(80,0)(60,0)
\dottedline{3}(80,0)(70,-10)(60,0)
\drawline[90](60,0)(50,20)(30,20)(20,30)
\dottedline{3}(20,30)(0,50)
\put(130,13){\scriptsize$p_2$}
}
\Thicklines
\dark{
\dottedline{3}(140,50)(120,30)
\drawline[90](120,30)(90,30)(80,10)(0,10)
\put(126,45){\scriptsize$p_1$}
}
\end{picture}
\end{center}
\caption{ }
\label{fig:p2}
\end{figure}

In all the cases either $\{\barred{n}, \barred{n-i}\} \in I$ or $\{n, \barred{n-i}\} \in I$ for some $i \geq 1$.
Note that if $s_{n-2}$ appears in between $s_n$ and $s_{n-1}$ in the expression of the Coxeter element $c$, that is, $c$ can be written as one of the following forms: $\cdots s_n \cdots s_{n-2} \cdots s_{n-1} \cdots$ or $\cdots s_{n-1} \cdots s_{n-2} \cdots s_{n} \cdots$, then the source of $p_1$ must be $n-1$. However, in this case, the indices $n-1$ and $\barred{n-1}$ form a single two cycle when $c$ is written as a permutation on the index set $[1, n] \cup [\barred{1}, \barred{n}]$ which implies that $\barred{n-1}$ does not belong to $c^m \cdot [1, k]$ for any $k \in [2, n-2]$ and $m \in \ZZ$.
On the other hand, if $s_{n-2}$ does not appear in between $s_n$ and $s_{n-1}$ in the expression of the Coxeter element $c$ then the indices $n$ and $\barred{n}$ form a single two cycle when $c$ is written as a permutation on the index set $[1, n] \cup [\barred{1}, \barred{n}]$. This implies that neither the index $n$ nor the index $\barred{n}$ belongs to $c^m \cdot [1, k]$ for any $k \in [2, n-2]$ and $m \in \ZZ$.
This shows that there does not exist a family of vertex-disjoint paths in $\Gamma (\mathsf{D_n}, c)$ with the sources and sinks labeled by~$c^m \cdot [1, k]$ for $k \in [2, n-2]$ such that there is exactly one crossing among its paths. This completes the proof of part $(1)$ of Proposition~\ref{pr:genminorsD}. In fact, it can be shown in a similar argument that no crossing (one or two) can appear in any family of vertex-disjoint paths in $\Gamma (\mathsf{D_n}, c)$ with the sources and sinks labeled by~$c^m \cdot [1, k]$ for $k \in [2, n-2]$.

Part $(2)$ and part $(3)$ of Proposition~\ref{pr:genminorsD} will become clear after we recall the corresponding spin representations.
Let $T$ be an $n$-subset of $[1, n] \cup [\barred{1}, \barred{n}]$, then the spin representations $V_{\omega_{n-1}}$ and $V_{\omega_{n}}$ can be realized as the vector spaces span by basis vectors as follows:
\begin{equation}
\label{eq:spinD}
\begin{aligned}
V_{\omega_{n-1}} &= \left\langle T \left| \begin{array}{lll}
        \mbox{ $i$ and $\overline{i}$ do not appear simultaneously in $T$, }\\
        \mbox{ there are an odd number of $\overline{i}$\,'s appearing in $T$. }
        \end{array} \right. \right\rangle \,,\\
V_{\omega_{n}} &= \left\langle T \left| \begin{array}{lll}
        \mbox{ with $i$ and $\overline{i}$ do not appear simultaneously in $T$, }\\
        \mbox{ there are an even number of $\overline{i}$\,'s appearing in $T$. }
        \end{array} \right. \right\rangle \,.\\
\end{aligned}
\end{equation}
The $\soe$-actions on $V_{\omega_{n-1}}$ and $V_{\omega_{n}}$ are given as follows:
\beal \label{soeactiononspin}
e_i \cdot T =
\bmcase
    T \setminus \{\,i+1, \overline{i}\,\} \cup \{\,i, \overline{i+1}\,\}, & \mbox{ if } i \neq n \mbox{ and } i+1, \, \overline{i} \, \in T\,; \\[0.1in]
    T \setminus \{\,\overline{n}, \overline{n-1}\,\} \cup \{n-1, n\}, & \mbox{ if } i = n \mbox{ and } \overline{n}, \, \overline{n-1} \, \in T\,; \\[0.04in]
    0, & \mbox{ otherwise} \ , \\
\emcase \\[.4in]
f_i \cdot T =
\bmcase
    T \setminus \{\,i, \overline{i+1}\,\} \cup \{\,i+1, \overline{i}\,\}, & \mbox{ if } i \neq n \mbox{ and } i, \, \overline{i+1} \, \in T\,; \\[0.1in]
    T \setminus \{n-1, n\} \cup \{\,\overline{n}, \overline{n-1}\,\}, & \mbox{ if } i = n \mbox{ and } n-1, \, n \, \in T\,; \\[0.04in]
    0, & \mbox{ otherwise} \ .\\
\emcase
\eea
Hence $x_i(t)$ and $x_\barred{i}(t)$ act as $I + t e_i$ and $I + t f_i$ on $V_{\omega_{n-1}}$ and $V_{\omega_n}$ respectively.
The fundamental representations $V_{\omega_{n-1}}$ and $V_{\omega_{n}}$ have highest weight vectors $\{1, 2, \ldots,
n-1, \overline{n}\}$ and $\{1, 2, \ldots, n-1, n\}$ respectively.

The combinatorial meaning of the graph $\Gamma (\mathsf{D_n}, c)$ in these cases is completely analogous to the one before.
The requirement of the paths being bundled is due to that the non-trivial actions of $e_i$ and $f_i$ require two specified indices to appear simultaneously in $T$. In this case, the coefficient of the corresponding basis vector should be $t$ instead of $t^2$, therefore we take the square root of the weight of a family of bundled vertex-disjoint paths. This completes the proof of Proposition~\ref{pr:genminorsD}.
\endproof

The proofs of the Propositions~\ref{pr:genminorsB}, \ref{pr:genminorsC} are similar to the proof of Proposition~\ref{pr:genminorsD}. 

\textbf{Proof of Proposition~\ref{pr:genminorsB}:} Let $\lg$ be the simple Lie algebra of type $\mathsf{B_n}$ for $n \geq 2$, that is, the odd special orthogonal Lie algebra $\so_{2n+1}$.  Then the action of generators in the standard $2n+1$-dimensional representation $V$ with respect to the standard basis $\vv_1, \ldots, \vv_n, \vv_0, \vv_{\barred{n}}, \ldots, \vv_{\barred{1}}$ can be written as:
\beal \label{sooaction}
e_i \cdot \vv_j =
\bmcase
    \vv_i\,, & \mbox{ if } i\neq n \mbox{ and } j = i+1\,;\\
    \vv_{\barred{i+1}}\,, & \mbox{ if } i \neq n \mbox{ and } j = \barred{i}\,;\\[.02in]
    \sqrt{2}\, \vv_n\,, & \mbox{ if } i = n \mbox{ and } j = 0\,;\\
    \sqrt{2}\, \vv_{0}\,, & \mbox{ if } i = n \mbox{ and } j = \barred{n}\,;\\
    0\,, & \mbox{ otherwise}, \\
\emcase \\[.5in]
f_i \cdot \vv_j =
\bmcase
    \vv_{i+1}\,, & \mbox{ if } i \neq n \mbox{ and } j = i\,; \\
    \vv_{\barred{i}}\,, & \mbox{ if } i \neq n \mbox{ and } j = \barred{i+1}\,;\\[.02in]
    \sqrt{2}\, \vv_0\,, & \mbox{ if } i = n \mbox{ and } j = n\,;\\
    \sqrt{2}\, \vv_{\barred{n}}\,, & \mbox{ if } i = n \mbox{ and } j = 0\,;\\
    0\,, & \mbox{ otherwise}, \\
\emcase
\eea
for $i \in [1, n]$ and $j \in [1, n] \cup \{0\} \cup [\barred{1}, \barred{n}]$.

It is easy to see that $x_i(t)$ and $x_\barred{i}(t)$ act as $I + t e_i + \frac{t^2}{2} e_i^2$ and $I + t f_i + \frac{t^2}{2} f_i^2$ respectively on $V$. The fundamental representation $V_{\omega_k}$ for $k = 1, \ldots, n-1$ is realized as $\bigwedge^k V$ with highest weight vector $\vv_1 \wedge \cdots \wedge \vv_{k}$.

To prove part $(1)$ of Proposition~\ref{pr:genminorsB}, it is enough to show that there is no crossing among any family of vertex-disjoint paths in $\Gamma (\mathsf{B_n}, c)$ with the sources and sinks labeled by~$c^m \cdot [1, k]$ for $k \in [2, n-1]$. Note in $\Gamma(\mathsf{B_n}, c)$, the crossing can only happen on the horizontal level $0$ in $x_n$ or $x_{\barred{n}}$ (see Figure~\ref{fig:bnchipsnsp}). The index $0$ does not belong to any index set corresponding to an extremal weight, and the only diagonal edges connecting the horizontal level $0$ are within the chips $x_n$ and $x_{\barred{n}}$ themselves. Together with the fact that $x_n$ and $x_{\barred{n}}$ only appear once in $\Gamma (\mathsf{B_n}, c)$, we conclude that such a crossing can not happen. This completes the proof of part $(1)$ of Proposition~\ref{pr:genminorsB}.

To prove part $(2)$ of Proposition~\ref{pr:genminorsB}, we first recall the spin representation in this case.
Let $T$ be an $n$-subset of $[1, n] \cup [\barred{1}, \barred{n}]$. Then the spin representation can be realized as a vector space span by basis vectors as follows:
\begin{equation}
\label{eq:spinB}
V_{\omega_{n}} = \left\langle T \left|
        \mbox{ $i$ and $\overline{i}$ do not appear simultaneously in $T$ }
        \right. \right\rangle \,.\\
\end{equation}
The $\soo$-action on $V_{\omega_{n}}$ is given as follows:
\beal \label{sooactiononspin}
e_i \cdot T =
\bmcase
    T \setminus \{\,i+1, \overline{i}\,\} \cup \{\,i, \overline{i+1}\,\}, & \mbox{ if } i \neq n \mbox{ and } i+1, \, \overline{i} \, \in T\,; \\[0.1in]
    T \setminus \{\,\overline{n}\,\} \cup \{n\}, & \mbox{ if } i = n \mbox{ and } \overline{n}\, \in T\,; \\[0.04in]
    0, & \mbox{ otherwise} \ , \\
\emcase \\[.4in]
f_i \cdot T =
\bmcase
    T \setminus \{\,i, \overline{i+1}\,\} \cup \{\,i+1, \overline{i}\,\}, & \mbox{ if } i \neq n \mbox{ and } i, \, \overline{i+1} \, \in T\,; \\[0.1in]
    T \setminus \{n\} \cup \{\,\overline{n}\,\}, & \mbox{ if } i = n \mbox{ and } n \, \in T\,; \\[0.04in]
    0, & \mbox{ otherwise} \ .\\
\emcase
\eea
Hence $x_i(t)$ and $x_\barred{i}(t)$ act as $I + t e_i$ and $I + t f_i$ on $V_{\omega_n}$ respectively.
The fundamental representation $V_{\omega_{n}}$ has highest weight vector $\{1, 2, \ldots, n-1, n\}$. The reasons for requiring the bundled condition and taking the square root of the weight of a collection of vertex-disjoint paths are as the same as those in the part $(2)$ and $(3)$ of Proposition~\ref{pr:genminorsD}.
This completes the proof of Proposition~\ref{pr:genminorsB}.
\endproof

\textbf{Proof of Proposition~\ref{pr:genminorsC}:} Let $\lg$ be the simple Lie algebra of type $\mathsf{C_n}$ for $n \geq 2$, that is, the symplectic Lie algebra $\sp_{2n}$.  Then the action of generators in the standard $2n$-dimensional representation $V$ with respect to the standard basis $\vv_1, \ldots, \vv_n, \vv_{\barred{n}}, \ldots, \vv_{\barred{1}}$ can be written as:
\beal \label{spaction}
e_i \cdot \vv_j =
\bmcase
    \vv_i\,, & \mbox{ if } i\neq n \mbox{ and } j = i+1\,;\\
    \vv_{\barred{i+1}}\,, & \mbox{ if } i \neq n \mbox{ and } j = \barred{i}\,;\\
    \vv_n\,, & \mbox{ if } i = n \mbox{ and } j = \barred{n}\,;\\
    0\,, & \mbox{ otherwise}, \\
\emcase \\[.5in]
f_i \cdot \vv_j =
\bmcase
    \vv_{i+1}\,, & \mbox{ if } i \neq n \mbox{ and } j = i\,; \\
    \vv_{\barred{i}}\,, & \mbox{ if } i \neq n \mbox{ and } j = \barred{i+1}\,;\\
    \vv_{\barred{n}}\,, & \mbox{ if } i = n \mbox{ and } j = n\,;\\
    0\,, & \mbox{ otherwise}, \\
\emcase
\eea
for $i \in [1, n]$ and $j \in [1, n] \cup [\barred{1}, \barred{n}]$.

It is easy to see that $x_i(t)$ and $x_\barred{i}(t)$ act as $I + t e_i$ and $I + t f_i$ respectively on $V$. Although the fundamental representation $V_{\omega_k}$ is not isomorphic to the exterior power $\bigwedge^k V$ for $k > 1$, it can be realized as a subrepresentation in $\bigwedge^k V$ with highest weight vector $\vv_1 \wedge \cdots \wedge \vv_{k}$. Hence, for our purpose, it makes no difference to work inside $\bigwedge^k V$. Proposition~\ref{pr:genminorsC} clearly holds since there is no crossing in any family of vertex-disjoint paths in $\Gamma(C_n, c)$. This completes the proof.
\endproof

\section*{Acknowledgments}
The author would like to thank A.~Zelevinsky for helpful comments and suggestions.


\end{document}